\RequirePackage{fix-cm}
\documentclass[smallextended,envcountsect]{svjour3}       
\smartqed  
\usepackage{graphicx}
\usepackage{amsmath}
\usepackage{subfig}
\usepackage{algorithm}
\usepackage{algorithmicx}
\usepackage{algpseudocode}
\usepackage{hyperref}
\usepackage{cleveref}
\crefname{algorithm}{Algorithm}{algorithms}
\crefname{chapter}{Chapter}{Chapters}
\crefname{equation}{}{Equations}
\crefname{figure}{Fig.}{Figures}
\crefname{lemma}{Lemma}{Lemmas}
\crefname{remark}{Remark}{Remarks}
\crefname{section}{Sect.}{Sections}
\crefname{subsection}{Sect.}{Sections}
\crefname{theorem}{Theorem}{theorems}

\makeatletter
\def\myunderbrace#1{\mathop{\vtop{\m@th\ialign{##\crcr
   $\hfil\displaystyle{#1}\hfil$\crcr
   \noalign{\kern3\p@\nointerlineskip}%
   \footnotesize\downbracefill\crcr\noalign{\kern3\p@\nointerlineskip}}}}\limits}
\makeatother

\graphicspath{{figures/}}

\newcommand{\bg}{\boldsymbol{g}}
\newcommand{\bu}{\boldsymbol{u}}
\newcommand{\bv}{\boldsymbol{v}}
\newcommand{\eps}{\epsilon}
\newcommand{\md}{\mathrm{d}}

\newcommand{\mB}{\mathcal{B}}
\newcommand{\mD}{\mathcal{D}}
\newcommand{\mF}{\mathcal{F}}

\newcommand{\mI}{\mathcal{I}}
\newcommand{\mL}{\mathcal{L}}
\newcommand{\mM}{\mathcal{M}}
\newcommand{\mN}{\mathcal{N}}
\newcommand{\mS}{\mathcal{S}}
\newcommand{\mP}{\mathcal{P}}
\newcommand{\mT}{\mathcal{T}}
\newcommand{\mO}{\mathcal{O}}
\newcommand{\lnorm}{\left \|}
\newcommand{\rnorm}{\right \|}
\newcommand{\diag}{\mathop{\mathrm{diag}}}
\newcommand{\udg}{\underline {g}}
\newcommand{\udA}{\underline {A}}
\newcommand{\udB}{\underline {B}}
\newcommand{\udU}{\underline {U}}
\newcommand{\veps}{\varepsilon}
\newcommand{\whu}{\widehat{u}}
\begin{document}

\title{Solving time-dependent PDEs with the ultraspherical spectral method
}


\author{Lu Cheng         \and
        Kuan Xu 
}


\institute{Lu Cheng \and Kuan Xu\at
            School of Mathematical Sciences, University of Science and Technology of China, 96 Jinzhai Road, Hefei 230026, Anhui, China \\
              \email{kuanxu@ustc.edu.cn}           
}

\date{ }

\maketitle

\begin{abstract}
  We apply the ultraspherical spectral method to solving time-\\dependent PDEs by proposing two approaches to discretization based on the method of lines and show that these approaches produce approximately same results. We analyze the stability, the error, and the computational cost of the proposed method. In addition, we show how adaptivity can be incorporated to offer adequate spatial resolution efficiently. Both linear and nonlinear problems are considered. We also explore time integration using exponential integrators with the 
  ultraspherical spatial discretization. Comparisons with the Chebyshev pseudospectral method are given along the discussion and they show that the ultraspherical spectral method is a competitive candidate for the spatial discretization of time-dependent PDEs.
\end{abstract}

\section{Introduction} \label{sec:intro}
In this article, we consider the one-dimension time-dependent PDE
\begin{subequations}
\begin{align}
\mT u &= \mF(t , u(x , t)), \label{tuF} \\
\mbox{s.t.}~~~ \mB u &= c, \label{Buc}  \\
u(x,0) &= f(x),
\end{align}
\label{tdp}%
\end{subequations}
where $\mT$ is the first-order differential operator in time. $\mF$ is a spatial operator that acts on the time $t$ and the solution $u(x, t)$. For a fixed $t$, $u(x, t)$ becomes a univariate function of the spatial variable $x$ defined on $[-1, 1]$. $\mF(t , u(x , t))$ can be further decomposed as
\begin{align}
\mF(t , u(x, t)) = \mL u + \mN(t, u(x, t)), \label{LN}
\end{align}
where $\mL$ and $\mN$ are the linear and nonlinear parts, respectively. Throughout this article, we follow the convention of writing $\mL u$, instead of $\mL (u)$, for $\mL$ is linear. Without loss of generality, we assume that $\mL$ is an $N$th order linear differential operator in space for $x \in [-1, 1]$
\begin{align}
\mL = a^N(x)\frac{{\md}^N}{{\md}x^N} + \ldots +a^1(x)\frac{{\md}}{{\md}x} + a^0(x) \label{opL}
\end{align}
with $a^N(x) \neq 0$ so that $\mL$ is non-singular. The side conditions $\mB$ contains $N$ linear functionals which are boundary conditions or constraints of other sorts and $c$ is an $N$-vector. The function $f(x)$ gives the initial condition.

In \cite{olv} Olver and Townsend present a fast and stable spectral method enabled by the ultraspherical polynomials which solves linear ordinary differential equations of the form
\begin{subequations}
\begin{align}
\mL u &= g, \label{luf}\\
\mbox{s.t. } \mB  u &= c , \label{buc}
\end{align}
\label{ode}%
\end{subequations}
where $\mL$ is also defined as in \eqref{opL}. This ultraspherical spectral method assumes the solution is written in its Chebyshev expansion
\begin{align*}
u(x) = \sum_{k=0}^{\infty}u_k T_k(x),
\end{align*}
where $T_k(x)$ is the Chebyshev polynomial of degree $k$. This way, $u(x)$ is identified by the coefficient vector $\bu = [u_0, u_1, \ldots]^{T}$. With a change of basis, the $\lambda$-order differentiation operator is as sparse as
\begin{align}
    \mD_{\lambda} = 2^{\lambda - 1}(\lambda - 1)!
    \begin{aligned}
      & ~~~\, \myunderbrace{\text{\footnotesize $\lambda$ times}} \\[-7pt]
      & \begin{pmatrix}
           \  0  \cdots  0  &~\lambda & & & &  \\
          & &\lambda +1&  & & \\
          & & &\lambda +2 & & \\
          & & & &\ddots& \\
      \end{pmatrix}
    \end{aligned},\label{opDiff}
\end{align}
for $\lambda = 1, 2, \ldots$, where $\mD_{\lambda}$ maps Chebyshev coefficients 
to ultraspherical $C^{(\lambda)}$ coefficients\footnote{In \cite{olv}, $\mD_0 = 
\mD_1$, while in this paper we let $\mD_1$ maps from Chebyshev $T$ to $C^{1}$ 
and $\mD_0 = \mI$, i.e., the identity operator, for notational consistency.}. 

If any of $a^{\lambda}(x)$ in \cref{opL} is not constant and written as 
\begin{align*}
a^{\lambda}(x) = \sum_{k=0}^{\infty}a_j C^{(\lambda)}_j(x), \nonumber
\end{align*}
the differential operator $\mD_{\lambda}$ should be pre-multiplied by the multiplication operator whose $(j,k)$ entry reads
\begin{align}
\mM_{\lambda}[a^{\lambda}]_{j, k}=\sum_{s=\max (0, k-j)}^{k} a_{2 s+j-k} c_{s}^{\lambda}(k, 2 s+j-k) \label{opMult}
\end{align}
for $j, k \geq 0$, where
\begin{align*}
c_{s}^{\lambda}(j, k)=\frac{j+k+\lambda-2 s}{j+k+\lambda-s} \frac{(\lambda)_{s}(\lambda)_{j-s}(\lambda)_{k-s}}{s !(j-s) !(k-s) !} \frac{(2 \lambda)_{j+k-s}}{(\lambda)_{j+k-s}} \frac{(j+k-2 s) !}{(2 \lambda)_{j+k-2 s}}.
\end{align*}
Note that $\mM_{\lambda}[a^{\lambda}]$ maps the ultraspherical space of $C^{(\lambda)}$ to itself. As long as $a^{\lambda}(x)$ possesses certain smoothness, it can be approximated by a finite series, that is,
\begin{align*}
a^{\lambda}(x) \approx \sum_{k=0}^{m}a_j C^{(\lambda)}_j(x). \nonumber
\end{align*}
This way, $\mM_{\lambda}[a^{\lambda}]$ becomes banded since $a_j = 0$ for $j > m$. Another approach to calculating the entries of $\mM_{\lambda}[a^{\lambda}]$ is given in \cite{tow1}, based on a recurrence relation for the multiplication operator.

When $\mD_{\lambda}$ and $\mM_{\lambda}[a^{\lambda}]$ are employed, each term in \cref{opL} maps to a different ultraspherical basis. So the following conversion operators $\mS_{\lambda}$ are needed to map the coefficients in $T$ to those in $C^{(1)}$ or $C^{(\lambda)}$ to $C^{(\lambda + 1)}$ respectively
\begin{subequations}
\begin{align}
\mS_0 &= \begin{pmatrix}
1&  &-\frac{1}{2}&  & & &\\
 &\frac{1}{2}& &-\frac{1}{2}& & & \\
 & &\frac{1}{2}& &-\frac{1}{2} & & \\
 & & &\ddots& &\ddots& \\
\end{pmatrix}, \label{S0} \\
\mS_{\lambda} &= \begin{pmatrix}
1&  &-\frac{\lambda}{\lambda +2}&  & & &\\
 &\frac{\lambda}{\lambda +1}& &-\frac{\lambda}{\lambda +3}& & & \\
 & &\frac{\lambda}{\lambda +2}& &-\frac{\lambda}{\lambda +4} & & \\
 & & &\ddots& &\ddots& \\
\end{pmatrix}
~\text{for }\lambda \geq 1.
\end{align}
\label{opConv}
\end{subequations}

In terms of \cref{opDiff}, \cref{opMult}, and \cref{opConv}, the differential equation \cref{luf} can be represented as
\begin{align}
\left( \mM_N[a^N] \mD_N + \sum_{\lambda = 0}^{N-1}{\mS}_{N-1}\ldots {\mS}_{\lambda}\mM_{\lambda}[a^{\lambda}]\mD_{\lambda} \right)\bu  = {\mS}_{N-1}\ldots {\mS}_0 \bg, \label{lufInf}
\end{align}
where $\bg$ is the vector containing the Chebyshev coefficients of $g(x)$. To make \cref{lufInf} of finite dimension, the operators are truncated by the projection operator $\mP_n = ({I}_n , \boldsymbol{0})$, where ${I}_n$ is the $n \times n$ identity matrix, with the dimension $n$ properly chosen. The truncated version of \cref{lufInf} reads
\begin{equation}
    \begin{split}
        \mP_{n-N} \left( \mM_N[a^N] \mD_N + \sum_{\lambda = 0}^{N-1}{\mS}_{N-1}\ldots {\mS}_{\lambda}\mM_{\lambda}[a^{\lambda}]\mD_{\lambda} \right)&\mP_{n}^{\top} \mP_{n}\bu  \\= \mP_{n-N}{\mS}_{N-1}&\ldots {\mS}_0\mP_{n}^{\top} \mP_{n}\bg,
    \end{split}\label{lufFin}
\end{equation}
where the unknown $\mP_{n}\bu$ and the (unconverted) right-hand side $\mP_{n}\bg$ are $n$-vectors and the differential operators on the left-hand side and the product of the conversion operators on the right-hand side are approximated by their truncated version of dimension $(n-N) \times n$ via \textit{exact truncation}. The system \cref{lufFin} is finally squared up to form an $n \times n$ system by the first $n$ columns of the discretized version of the boundary conditions \cref{buc} and this is the system by solving which one obtains the Chebyshev coefficients $u_k$ of the truncated version of the solution
\begin{align*}
\widetilde{u}_n(x) = \sum_{k=0}^{n-1}u_k T_k(x). 
\end{align*}

The ultraspherical spectral method recapitulated above enjoys a few important advantages over the collocation-based pseudospectral methods, including linear computational complexity, good conditioning, and adaptivity via optimal truncation.

In this article, we extend the ultraspherical spectral method to the solution of the time-dependent problem \cref{tdp} within the method of lines (MOL) framework. Our investigation is by no means the first attempt to solve time-dependent PDEs by the ultraspherical spectral method. In \cite{tow2}, 
Townsend and Olver describe an extension of the ultraspherical spectral method to two spatial dimensions for the solution of linear PDEs with variable coefficients defined on bounded rectangular domains and their focus is on the automated manner of solution provided that the splitting rank of the partial differential operator (PDO) is known. When applied to an initial boundary value problem, this bivariate ultraspherical spectral method treats it as a boundary 
value problem of two spatial dimensions by deeming the time variable as a second spatial variable. Our motivation in this article, however, is to employ the ultraspherical spectral method in space while do the time-stepping using common time integration schemes. Moreover, we consider a more general setting where the problem may or may not have a sufficiently concise closed-form 
description or the spatial operator can only be evaluated via black-box routines, which is often the case in real-world problems. 

The first and probably only existing works where the ultraspherical spectral method is used in conjunction with time-stepping schemes may be \cite{slo,for}, where the implicit-explicit method and the backward Euler method are 
employed, respectively. However, the application of these time-stepping schemes are not theoretically analyzed to give insights on their performance. On the software side, the \textsc{Dedalus} package solves time-dependent PDEs using (a first-order variant of) the ultraspherical spectral method with the time-integration done by a range of ODE integrators including multistep and Runge-Kutta IMEX methods \cite{bur}. The success of these attempts suggests a pressing demand on the theoretical analysis of time stepping when the ultraspherical spectral method is used for the spatial discretization. This is exactly what the present article focuses on. By giving a rather complete treatment to solving time-dependent PDEs in one spatial dimension, this article may well serve as a foundation for migration to problems in higher spatial dimensions.

In the first part of this article, we concentrate on the linear case of 
\eqref{tdp}, i.e.,
\begin{subequations}
\begin{align}
\mT u &= \mL u, \label{tulu} \\
\mbox{s.t.}~~~ \mB u &= c, \label{Bucl}  \\
u(x,0) &= f(x),
\end{align}
\label{tdpl}%
\end{subequations}
by discussing the discretization of \cref{tdpl} via standard time stepping schemes (\cref{sec:disc}) and analyzing the stability (\cref{sec:stability}), the error (\cref{sec:error}), and the computational cost (\cref{sec:cost}). The stepping nature of the method enables an adaptive implementation which we describe in \cref{sec:adapt}. In the linear regime, our discussion will make frequent use of the one-dimensional transport equation 
\begin{subequations}
\begin{align}
u_t(x , t) &= u_x(x , t),  \\
u(1 , t) &= 0  
\end{align}
\label{wave1}
\end{subequations}
and the heat equation
\begin{subequations}
\begin{align}
u_t(x, t) &= u_{xx}(x, t),  \\
u(-1, t) &=u(1, t) = 0, 
\end{align}
\label{heat}%
\end{subequations}
both subject to the initial condition $u(x, 0) = f(x)$. Also, we shall simply 
take $f(x) = \exp(-200x^2)$ and $f(x) = \sin(2\pi x)$ for \cref{wave1} and 
\cref{heat}, respectively. In the study of the collocation-based pseudospectral 
method, much attention has been paid to these problems from various 
perspectives, particularly regarding the  stability restrictions on time 
stepping and the eigenvalue distribution of the spatial discretization 
operators, see, e.g., \cite{dub,got,wei,sne}.

We close our discussion in the linear regime by briefly analyzing the problems 
with periodic boundary conditions (\cref{sec:periodic}). The collocation-based 
pseudospectral method, for many years, has been taken as `the' method, and the 
discussion and analysis for the linear case facilitate the comparison between 
the two methods. In addition, they lay the foundation for the analysis of 
nonlinear time-dependent problems (\cref{sec:nonlinearity}). In \cref{sec:ei}, 
we examine the application of the exponential integrator in conjunction with 
the ultraspherical spectral method. Conclusion and discussion are given in the 
final section.

Throughout this article, all the norms are taken to be the infinity norm. 
Calligraphy font is used for operators or infinite matrices and bold fonts for 
infinite vectors, whereas the truncated version of operators, infinite 
matrices, and vectors are in normal fonts.

All the numerical experiments in this article are performed in \textsc{Julia} 
v1.5.3 on a desktop with a 4 core 2.1 Ghz AMD Ryzen 5 3500U CPU.

\section{Discretization} \label{sec:disc}
We start by considering the discretization of \cref{tdpl}. Suppose that the 
solution $u(x, t)$ is written as an infinite Chebyshev series
\begin{align*}
u(x, t) = \sum_{k=0}^{\infty}u_k(t)T_k(x), 
\end{align*}
where the coefficients $u_k(t)$ we are solving for are now dependent of time $t$. If the spatial operator $\mL$ on the right-hand side of \cref{tulu} is expressed in terms of the operators reviewed in \cref{sec:intro} as
\begin{align}
\mL=\mM_{N}\left[a^N\right] \mD_{N}+\sum_{\lambda=0}^{N-1} \mS_{N-1} \cdots \mS_{\lambda} \mM_{\lambda}\left[a^{\lambda}\right] \mD_{\lambda}, \label{L}
\end{align}
the left-hand side of \cref{tulu} must be pre-multiplied by a series of conversion operators so that both the sides end up being in the $C^{(N)}$ basis, that is,

\begin{align}
{\mS}_{N-1}\ldots \mS_0\mT \bu = 
\mL \bu, \label{tulu_us}
\end{align}
where $\bu = [u_0(t), u_1(t), \ldots]^{T}$ is the infinite vector collecting the coefficients $u_k(t)$.

Now we bifurcate our discussion by presenting two ways to further discretize 
\cref{tulu_us} and enforce the boundary condition \cref{Bucl}, both following 
the method of lines. They differ in how a square system is formed by solving 
which we obtain a truncated approximation to $\bu$.

In the remainder of this article, we confine our discussion about the discretization of the temporal operator $\mT$ to the standard time marching schemes for solving the ODE initial value problem $v_t = f(t,v)$. That is, we consider the linear multistep methods 
\begin{align}
\sum_{j = 0}^{r}\alpha_j {v}^{k + j} = h\sum_{j = 0}^{r}\beta_j f^{k + j}, \label{LMM}
\end{align}
where $\alpha_r=1$, and the explicit Runge-Kutta methods
\begin{subequations}
\begin{align}
y_j &= h f(t_k + \theta_j h , v^k + \mu_j y_{j-1}), \mbox{ for } j = 1 , 2 , \cdots , s \label{os1}\\
v^{k+1} &= v^{k} + \sum_{j = 1}^{s}\gamma_j y_j, 
\end{align} \label{RK}%
\end{subequations}
where $\theta_1 = \mu_1 = 0$ and $\sum\limits_{j=1}^s \gamma_j = 1$. In \cref{LMM} and \cref{RK}, $h$ is the step size.

\subsection{Approach 1} Our first approach enforces the main equation and the boundary conditions simultaneously. To this end, we truncate the operators and $\bu$
\begin{align}
\mP_{n-N}{\mS}_{N-1}\ldots \mS_0 \mP_{n}^{\top}\mT  \mP_{n}\bu = 
\mP_{n-N} \mL \mP_{n}^{\top}\mP_{n} \bu, \label{tulu_truncated}
\end{align}
which amounts to taking the first $n-N$ rows and the first $n$ columns of ${\mS}_{N-1}\ldots \mS_0$ and $\mL$ and approximating the solution by its $n$-term truncation
\begin{align*}
u_n(x, t) \approx \widetilde{u}_n(x, t) = \sum_{k=0}^{n-1}u_k(t)T_k(x).
\end{align*}
Note that the truncation of $\mP_{n-N} \mL \mP_{n}^{\top}$ is done \textit{exactly} as
\begin{align*}
\mP_{n-N} \mL \mP_{n}^{\top} =& \mP_{n-N} \left( \mM_N[a^N] \mD_N + \sum_{\lambda = 0}^{N-1}{\mS}_{N-1}\ldots {\mS}_{\lambda}\mM_{\lambda}[a^{\lambda}]\mD_{\lambda} \right)\mP_{n}^{\top} \\ 
=& \left( \mP_{n-N} \mM_N[a^N] \mP_{n-N}^{\top}\right) \left( \mP_{n-N} \mD_N \mP_{n}\right) + \sum_{\lambda = 0}^{N-1}\left( \mP_{n-N}{\mS}_{N-1}\mP_{n - N + 2}^{\top} \right) 
\\
\times & \left( \prod_{i=2}^{N - \lambda} \mP_{n - N+ 2(i-1)}{\mS}_{N-i}\mP_{n - N  + 2i}^{\top}\right) \left(\mP_{n - N + 2(N - \lambda)}\mM_{\lambda}[a^{\lambda}] \mP_{n-\lambda}^{\top}\right) \\
\times & \left(\mP_{n-\lambda}\mD_{\lambda} \mP_{n}^{\top}\right).
\end{align*}
For exact truncations of operators, see \cite[Remark 2]{olv} for details.

\begin{figure}[tbhp]
    \centering
    \subfloat[Explicit linear multistep methods]{\label{fig:disc1_a}\includegraphics[scale=0.75]{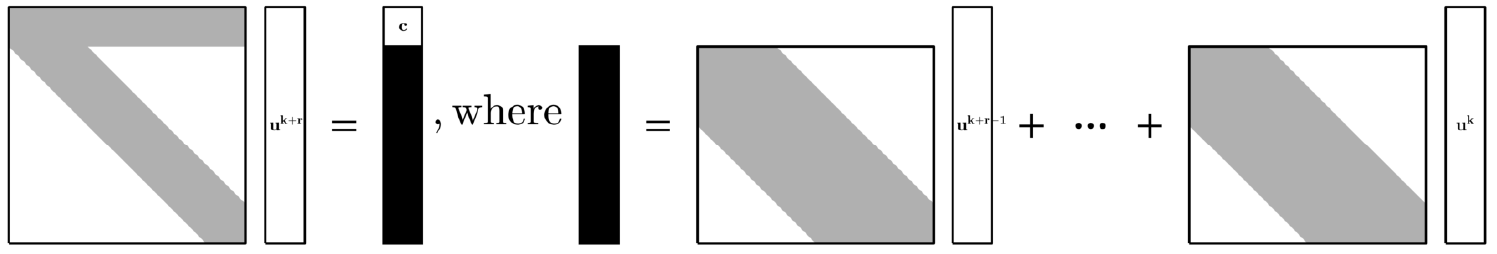}}
    
    \subfloat[Implicit linear multistep methods]{\label{fig:disc1_b}\includegraphics[scale=0.75]{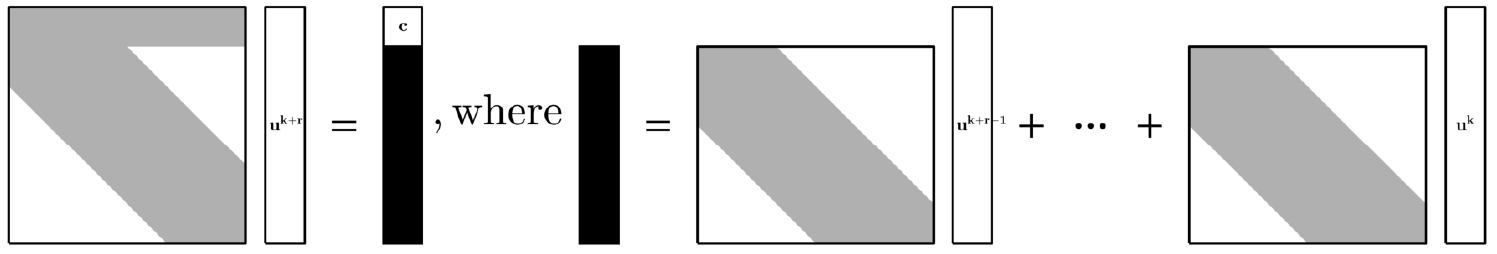}}
    
    \subfloat[Runge-Kutta methods]{\label{fig:disc1_c}\includegraphics[scale=0.75]{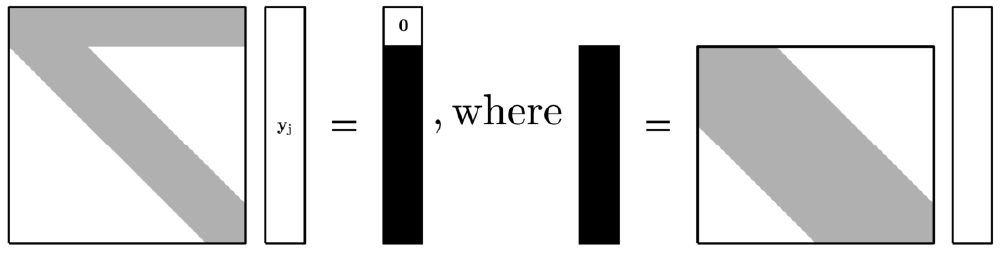}}
    
    \caption{Sparsity patterns of the fully discretized systems in Approach 1 for 
    linear multistep methods \cref{disc1_lms} and Runge-Kutta methods 
    \cref{disc1_rk}. \cref{fig:disc1_a} and \cref{fig:disc1_b} mainly differ in the 
    matrix on the left-hand side in that the lower bandwidth of its banded part is 
    zero for explicit schemes whereas the banded part could have nonzero 
    sub-diagonals for implicit schemes. }
    \label{fig:disc1}
\end{figure}

We truncate the operators in the boundary conditions analogously by taking the 
first $n$ columns of $\mB$
\begin{align}
\mB\mP_{n}^{\top}\mP_{n}\bu = c. \label{bc_truncated}
\end{align}
When \cref{bc_truncated} is laid on the top of \cref{tulu_truncated}, an $n \times n$ square system is formed despite that the temporal operator is not yet discretized.

Now we turn to the discretization in time. When a multistep method is applied to \cref{tulu_truncated}, we have
\begin{align*}
\sum_{j = 0}^{r}\alpha_j \mP_{n-N}{\mS}_{N-1}\ldots \mS_0 \mP_{n}^{\top} \mP_{n}\bu^{k + j} = h\sum_{j = 0}^{r}\beta_j \mP_{n-N} \mL \mP_{n}^{\top}\mP_{n} \bu^{k + j}.
\end{align*}
or, equivalently,
\begin{equation}
    \begin{split}
        (& \mP_{n-N}{\mS}_{N-1}\ldots \mS_0 \mP_{n}^{\top} - h\beta_r\mP_{n-N} \mL \mP_{n}^{\top})\mP_{n}\bu^{k + r} \\ &= h\sum_{j = 0}^{r - 1}\beta_j \mP_{n-N} \mL \mP_{n}^{\top}\mP_{n} \bu^{k + j} - \sum_{j = 0}^{r-1}\alpha_j \mP_{n-N}{\mS}_{N-1}\ldots \mS_0 \mP_{n}^{\top} \mP_{n}\bu^{k + j}. 
    \end{split} \label{tulu_lms}
\end{equation}
Here, ${\bu}^{k} = [u_0(t_k) , u_1(t_k) , \ldots ]^{T}$ is the approximate solution at $k$th time step, and $\mP_{n}\bu^{k}$ is the $n$-vector with the trailing coefficients dropped.

When a Runge-Kutta method is applied to \cref{tulu_truncated}, each stage becomes
\begin{align}
\mP_{n-N}{\mS}_{N-1}\ldots \mS_0 \mP_{n}^{\top} y_j = h \mP_{n-N}\mL \mP_{n}^{\top}(\mP_{n}\bu^{k} + \mu_jy_{j-1}) \label{tulu_rk}
\end{align}
for $j = 1 , 2 , \cdots , s$. Note that $y_j$ is a finite vector (see \eqref{os1}).

We are finally in a position to form the fully discretized square system. Stacking the boundary conditions \cref{bc_truncated} and the main equation \cref{tulu_lms} gives
\begin{equation}
    \begin{split}
        & \begin{pmatrix}
        \mB\mP_{n}^{\top}\\ 
        \mP_{n-N}{\mS}_{N-1}\ldots \mS_0 \mP_{n}^{\top} - h\beta_r\mP_{n-N} \mL \mP_{n}
        \end{pmatrix} \mP_{n}\bu^{k+r} \\
        & = \begin{pmatrix}
        c \\ 
        \sum\limits_{j = 0}^{r - 1}(\beta_j h\mP_{n-N} \mL\mP_{n}^{\top} - \alpha_j \mP_{n-N}{\mS}_{N-1}\ldots \mS_0 \mP_{n}^{\top}) \mP_{n}\bu^{k + j}
        \end{pmatrix}, 
    \end{split}\label{disc1_lms}
\end{equation}
by solving which we have $\mP_{n}\bu^{k+r}$. The sparsity structure of \cref{disc1_lms} is shown by \cref{fig:disc1_a} and \cref{fig:disc1_b} for explicit and implicit multistep methods, respectively.

When \cref{tulu_rk} is combined with the boundary conditions, we obtained a square system for the intermediate solutions $y_j$ at each stage of a Runge-Kutta method
\begin{align}
\begin{pmatrix}
\mB\mP_{n}^{\top}\\
\mP_{n-N}{\mS}_{N-1}\ldots \mS_0 \mP_{n}^{\top}
\end{pmatrix}y_j &= \begin{pmatrix}
0 \\
h\mP_{n-N}\mL \mP_{n}^{\top}(\mP_{n}\bu^{k} + \mu_jy_{j-1})
\end{pmatrix},
\label{disc1_rk}
\end{align}
where $j = 1, 2, \cdots, s$, and we update $\mP_{n}^{\top} \bu^{k+1}$ as
\begin{align*}
\mP_{n} \bu^{k+1} &= \mP_{n} \bu^{k} + \sum_{j = 1}^{s}\gamma_jy_j.
\end{align*}
The sparsity of \cref{disc1_rk} is shown by \cref{fig:disc1_c}.

\subsection{Approach 2} 
\begin{figure}[tbhp]
    \centering
    \subfloat[Explicit linear multistep methods]{\label{fig:disc2_a}\includegraphics[scale=0.6]{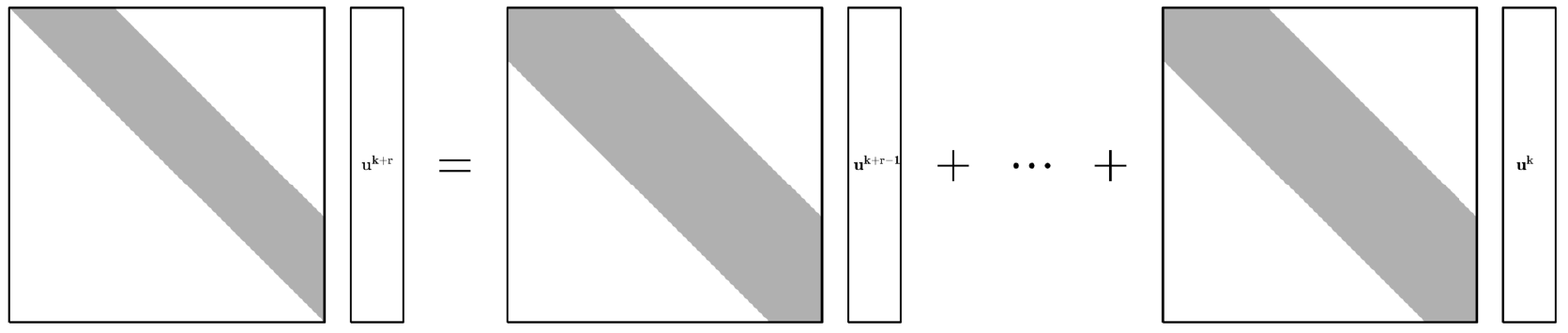}}
    
    \subfloat[Implicit linear multistep methods]{\label{fig:disc2_b}\includegraphics[scale=0.6]{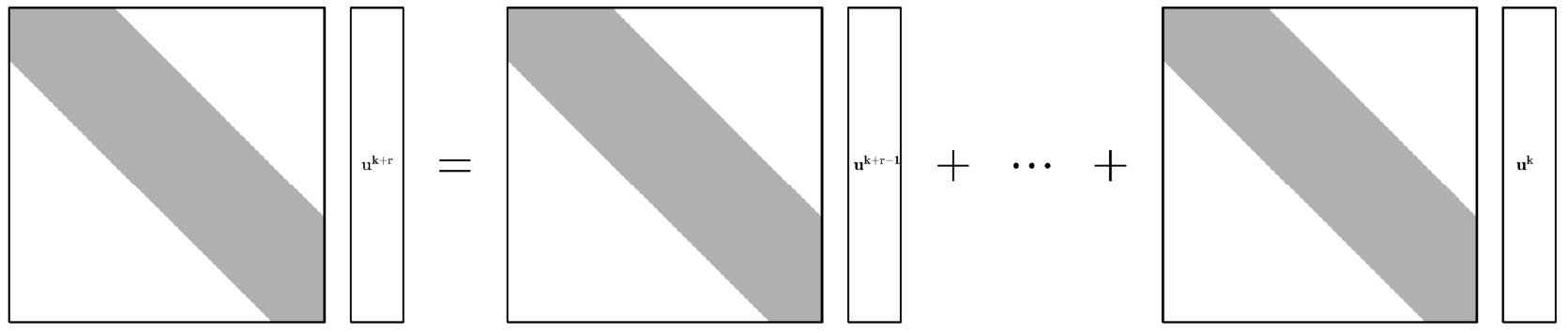}}
    
    \subfloat[Runge-Kutta methods]{\label{fig:disc2_c}\includegraphics[scale=0.364]{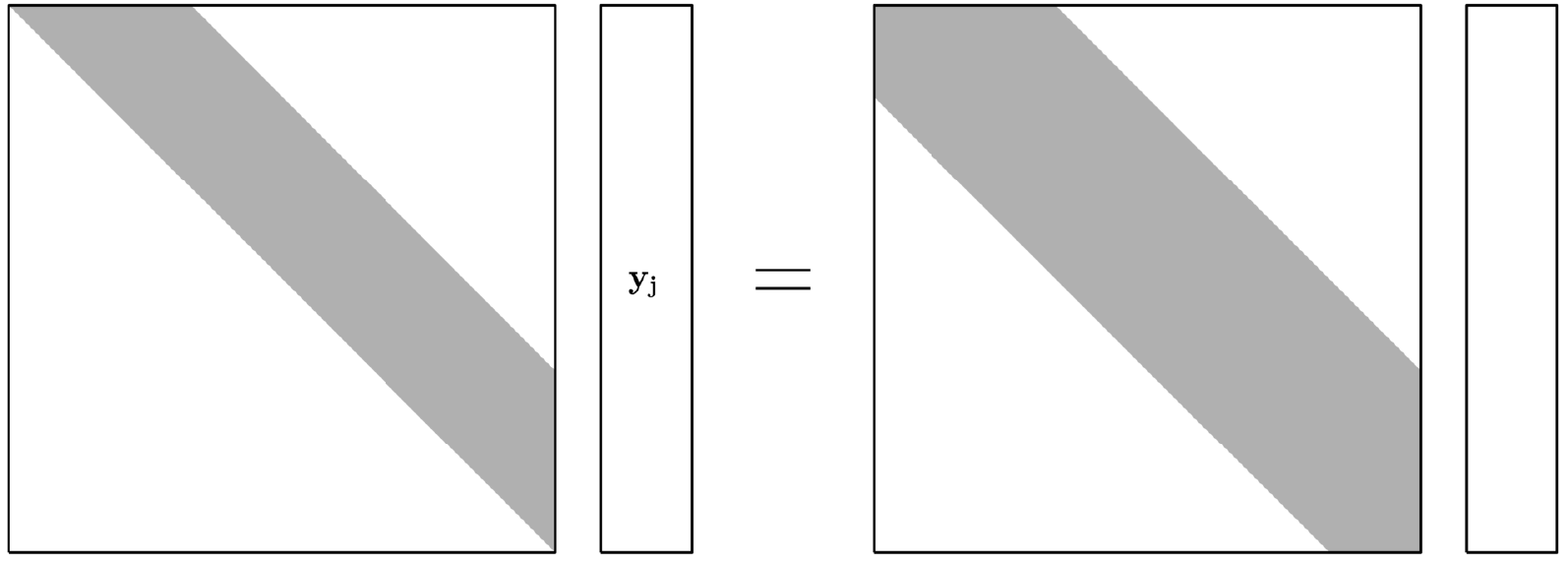}}
    
    \caption{Sparsity patterns of the fully discretized system in Approach 2 for 
    linear multistep methods \cref{tulu_lms_a2} and Runge-Kutta methods 
    \cref{tulu_rk_a2}.}
    \label{fig:disc2}
\end{figure}

Approach 2 ignores the boundary conditions in the first place and truncates \cref{tulu_us} to form a square system:
\begin{align}
\mP_{n}{\mS}_{N-1}\ldots \mS_0 \mP_{n}^{\top} \mT\mP_{n}\bu = \mP_{n} \mL \mP_{n}^{\top}\mP_{n} \bu,\label{tulu_truncated2}
\end{align}
where the truncations of $\mP_{n}{\mS}_{N-1}\ldots \mS_0 \mP_{n}^{\top}$ and $\mP_{n} \mL \mP_{n}^{\top}$ are, again, carried out exactly. Stepping using multistep or Runge-Kutta methods, we end up with
\begin{equation}
    \begin{split}
        (\mP_{n}&{\mS}_{N-1}\ldots \mS_0 \mP_{n}^{\top} - h\beta_r\mP_{n} \mL \mP_{n})\mP_{n}\bu^{k + r} \\
        & = h\sum_{j = 0}^{r - 1}\beta_j \mP_{n} \mL \mP_{n}^{\top}\mP_{n} \bu^{k + j} - \sum_{j = 0}^{r-1}\alpha_j \mP_{n}{\mS}_{N-1}\ldots \mS_0 \mP_{n}^{\top} \mP_{n}\bu^{k + j} 
    \end{split} \label{tulu_lms_a2}
\end{equation}
and
\begin{align}
\mP_{n}{\mS}_{N-1}\ldots \mS_0 \mP_{n}^{\top} y_j = h \mP_{n}\mL \mP_{n}^{\top}(\mP_{n}\bu^{k} + \mu_j y_{j-1}), \label{tulu_rk_a2}
\end{align}
respectively. The sparsity patterns are shown in \cref{fig:disc2}. Note that \cref{tulu_lms_a2} and \cref{tulu_rk_a2} differ from \cref{tulu_lms} and \cref{tulu_rk} by being square systems instead of rectangular. Since $\mP_{n}{\mS}_{N-1}\ldots \mS_0 \mP_{n}^{\top}$ is non-singular, \cref{tulu_lms_a2} and \cref{tulu_rk_a2} can be solved for $\mP_{n} \bu^{k+r}$ and the intermediate result at $j$th stage, respectively. Obviously, $\mP_{n} \bu^{k+r}$ obtained this way rarely satisfies the boundary condition. To enforce the boundary condition, we free $N$ components in $\mP_{n} \bu^{k+r}$ and allow them to be re-determined by
\begin{align*}
\mB\mP_{n}^{\top} \mP_{n}\bu^{k+r} = c.
\end{align*}
This is, in fact, an $N \times N$ system. For example, if we choose to 
re-determine the last $N$ components in $\mP_{n} \bu^{k+r}$, we end up with the 
$N \times N$ system in, for example, \textsc{Matlab}'s syntax
\begin{align*}
\mB\texttt{(1:N,n-N+1:n)}\bu^{k+r}\texttt{(n-N+1:n)} = 
c  - \mB\texttt{(1:N,1:n-N)}\bu^{k+r}\texttt{(1:n-N)}.
\end{align*}

One may wonder the difference between the solutions obtained via these two approaches, which we investigate now.

\subsection{Approach 1 vs Approach 2}\label{app}
Assuming the solution of \cref{tdpl} is Lipschitz continuous in both space and 
time, we now show that the difference between the solutions obtained via the 
two approaches is bounded and vanishes when the discretization in both time and 
space becomes increasingly dense. Although what is furnished below is not a 
rigorous proof, it suffices to give an explanation why these two approaches 
return converging solutions.

We set out by considering three problems and assuming the solutions $u^{k+j}$ are known for $j = 0, 1, \ldots, r-1$.
\begin{itemize}
\item Problem 1: use Approach 1 to obtain a solution vector of length $n+N$, which can be computed by solving the following system (see \cref{disc1_lms}, \cref{disc1_rk}, and \cref{fig:disc1}):
\begin{align}
\begin{pmatrix}
B_{(1)}& B_{(2)}\\
S_{(1)}& S_{(2)}\\
S_{(3)}& S_{(4)}
\end{pmatrix}u_{P_1}^{k+r}=
\begin{pmatrix}
c\\
\sum\limits_{j=0}^{r-1}\begin{pmatrix}
L^j_{(1)}& L^j_{(2)}\\
L^j_{(3)}& L^j_{(4)}
\end{pmatrix}u^{k+j}
\end{pmatrix}. \label{A1P1}
\end{align}

Note that \cref{A1P1} covers both linear multistep and Runge-Kutta methods. On the left-hand side, the top $N$ rows are partitioned as an $N \times n$ part $B_{(1)}$ and an $N \times N$ part $B_{(2)}$. Along with the $N$-vector $c$, the first $N$ equations represent the boundary conditions. The banded part of the coefficient matrix is partitioned into $S_{(1)}$, $S_{(2)}$, $S_{(3)}$, and $S_{(4)}$, whose dimensions are $(n-N) \times n$, $(n-N) \times N$, $N \times n$, and $N \times N$, respectively. The solution, denoted by $u^{k+r}_{P_1}$, is an $(n+N)$-vector. The $(n+N)$-vector $u^{k+j}$ represents either the solution at the previous steps or the initial condition. The banded matrices which $u^{k+j}$ multiplies with are partitioned into $L^j_{(1)}$, $L^j_{(2)}$, $L^j_{(3)}$, and $L^j_{(4)}$, whose dimensions are conformal with $S_{(1)}$, $S_{(2)}$, $S_{(3)}$, and $S_{(4)}$, respectively.
\item Problem 2: use Approach 1 to obtain a solution vector of length $n$, which amounts to solving the $n \times n$ system
\begin{align*}
\begin{pmatrix}
B_{(1)}\\
S_{(1)}
\end{pmatrix}u_{P_2}^{k+r}=
\begin{pmatrix}
c\\
\sum\limits_{j=0}^{r-1}L^j_{(1)}u^{k+j}\texttt{(1:n)}
\end{pmatrix}.
\end{align*}
\item Problem 3: use Approach 2 to obtain a solution vector of length $n$ by first solving the $n \times n$ system
\begin{subequations}
\begin{align}
\begin{pmatrix}
S_{(1)}\\
S_{(3)}
\end{pmatrix}\whu_{P_3}^{k+r}&=\sum\limits_{j=0}^{r-1}
\begin{pmatrix}
L^j_{(1)}\\
L^j_{(3)}
\end{pmatrix}u^{k+j}\texttt{(1:n)}, \label{P3}
\end{align}
which produces the intermediate solution $\whu_{P_3}^{k+r}$. This is then followed by the correction step which re-determines the last $N$ components of $\whu_{P_3}^{k+r}$. If the corrected solution is denoted by $u_{P_3}^{k+r}$, the boundary conditions are satisfied
\begin{align}
B_{(1)}u_{P_3}^{k+r}=c. 
\end{align}
\end{subequations}
\end{itemize}

Now we assume that $n$ is large enough so that the solution is fully resolved and spectral accuracy is achieved in space. Thus, there exists a small number $\epsilon > 0$ so that $\lnorm u^{k+j}\texttt{(n-N+1:n+N)}\rnorm < \epsilon$ for $j = 0, 1, \ldots, r-1$ and $\lnorm u_{P_1}^{k+r}\texttt{(n-N+1:n+N)}\rnorm < \epsilon$. Also, we assume $h$ is small enough to stabilize whatever time stepper we are using. 

In the following argument, Problem 1 serves as a bridge connecting Problems 2 and 3 whose solutions are what we try to show to be close. Thus, we bound the difference between $u_{P_1}^{k+r}$ and $u_{P_3}^{k+r}$ (Step 1) and that between $u_{P_1}^{k+r}$ and $u_{P_2}^{k+r}$ (Step 2) first, and these results, when combined, give the difference between $u_{P_2}^{k+r}$ and $u_{P_3}^{k+r}$.

\textbf{Step 1:}
We first look at the difference between the first $n$ components of $u_{P_1}^{k+r}$ and the uncorrected solution $\whu_{P_3}^{k+r}$, i.e., $e_1 = ||u_{P_1}^{k+r}\texttt{(1:n)}-\whu_{P_3}^{k+r}||$.
From \cref{A1P1}, we have
\begin{align*}
\begin{pmatrix}
S_{(1)}& S_{(2)}\\
S_{(3)}& S_{(4)}
\end{pmatrix}u_{P_1}^{k+r}=\sum\limits_{j=0}^{r-1}
\begin{pmatrix}
L^j_{(1)}& L^j_{(2)}\\
L^j_{(3)}& L^j_{(4)}
\end{pmatrix}u^{k+j},
\end{align*}
that is,
\begin{align*}
\begin{pmatrix}
S_{(1)}\\
S_{(3)}
\end{pmatrix}u_{P_1}^{k+r}\texttt{(1:n)} & +
\begin{pmatrix}
S_{(2)}\\
S_{(4)}
\end{pmatrix}u_{P_1}^{k+r}\texttt{(n+1:n+N)}\\
& =  \sum\limits_{j=0}^{r-1}
\begin{pmatrix}
L^j_{(1)}\\
L^j_{(3)}
\end{pmatrix}u^{k+j}\texttt{(1:n)}+ \sum\limits_{j=0}^{r-1}
\begin{pmatrix}
L^j_{(2)}\\
L^j_{(4)}
\end{pmatrix}u^{k+j}\texttt{(n+1:n+N)}.
\end{align*}
Combining the last equation with \cref{P3}, we have
\begin{align}
e_1&=\lnorm \begin{pmatrix}
S_{(1)}\\
S_{(3)}
\end{pmatrix}^{-1}
\begin{pmatrix}
S_{(2)}u_{P_1}^{k+r}\texttt{(n+1:n+N)} - \sum\limits_{j=0}^{r-1}L^j_{(2)}u^{k+j}\texttt{(n+1:n+N)}\\
S_{(4)}u_{P_1}^{k+r}\texttt{(n+1:n+N)} - \sum\limits_{j=0}^{r-1}L^j_{(4)}u^{k+j}\texttt{(n+1:n+N)}
\end{pmatrix}\rnorm \nonumber \\
&\le \lnorm \begin{pmatrix}
S_{(1)}\\
S_{(3)}
\end{pmatrix}^{-1} \rnorm \max_j\{\|S_{(2)}\|,\|L^j_{(2)}\|,\|S_{(4)}\|,\|L^j_{(4)}\| \} 2\epsilon = C_1 \epsilon. \label{e1}
\end{align}

Now we bound the correction due to the enforcement of the boundary conditions, i.e., $e_2 = \lnorm \whu_{P_3}^{k+r}-u_{P_3}^{k+r}\rnorm$. 

For a multistep method \cref{LMM}, we have 
\begin{align}
B_{(1)}\whu_{P_3}^{k+r}-B_{(1)}{u}_{P_3}^{k+r} = B_{(1)}\whu_{P_3}^{k+r}-c = &B_{(1)}\whu_{P_3}^{k+r} + B_{(1)}\sum_{j = 0}^{r-1}\frac{\alpha_j}{\alpha_r}u^{k+j}, \label{e2_pre}
\end{align}
where we have used the fact that all $u_{P_3}^{k+j}$ satisfy the boundary conditions, i.e., $B^{(1)}u_{P_3}^{k+j} = c$ for $j = 0, \ldots, r-1$, and the consistency condition
\begin{align*}
\sum_{j = 0}^{r}\alpha_j = 0.
\end{align*}
Substituting \cref{LMM} into \cref{e2_pre} gives
\begin{align*}
& B_{(1)}\whu_{P_3}^{k+r}-B_{(1)}{u}_{P_3}^{k+r} \\ 
& = 
\frac{B_{(1)}}{\alpha_r}h 
\left(\sum_{j = 0}^{r-1}\beta_j
\begin{pmatrix}
    S_{(1)}^{j}\\
    S_{(3)}^{j}
\end{pmatrix}^{-1}
\begin{pmatrix}
L_{(1)}^{j}\\
L_{(3)}^{j}
\end{pmatrix}u^{k+j} + 
\beta_r
\begin{pmatrix}
    S_{(1)}^{r}\\
    S_{(3)}^{r}
\end{pmatrix}^{-1}
\begin{pmatrix}
L_{(1)}^{r}\\
L_{(3)}^{r}
\end{pmatrix}\whu_{P_3}^{k+r}\right) ,
\end{align*}
which further leads to 
\begin{align}
e_2 = \lnorm\frac{h}{\alpha_r}\left(\sum_{j = 0}^{r-1}\beta_j
\begin{pmatrix}
L_{(1)}^{j}\\
L_{(3)}^{j}
\end{pmatrix}u^{k+j} + 
\beta_r\begin{pmatrix}
L_{(1)}^{r}\\
L_{(3)}^{r}
\end{pmatrix}\whu_{P_3}^{k+r}\right)\rnorm \le C_2h, \label{e2_lmm}
\end{align}

Analogously, for Runge-Kutta methods \cref{RK}
\begin{align*}
B_{(1)}\whu_{P_3}^{k+r}-B_{(1)}{u}_{P_3}^{k+r} = B_{(1)}\whu_{P_3}^{k+r}-c = B_{(1)}\whu_{P_3}^{k+r} - B_{(1)}u^{k+r-1} = B_{(1)}\sum_{j = 1}^{r}\gamma_jy_j,
\end{align*}
implying
\begin{align}
e_2 = \lnorm\sum_{j = 1}^{r}\gamma_jy_j\rnorm \le C_3 h. \label{e2_rk}
\end{align}

Finally, \cref{e1} and any one of \cref{e2_lmm} and \cref{e2_rk} give
\begin{align}
e_3 &= \lnorm \begin{pmatrix}
u_{P_3}^{k+r}\\
0
\end{pmatrix} - u_{P_1}^{k+r}\rnorm \nonumber \\
&\le \lnorm
\begin{pmatrix}
\widehat{u}_{P_3}^{k+r}\\
0
\end{pmatrix} - 
\begin{pmatrix}
u_{P_1}^{k+r}\texttt{(1:n)}\\
0
\end{pmatrix} \rnorm + \lnorm 
\begin{pmatrix}
u_{P_1}^{k+r}\texttt{(1:n)}\\
0
\end{pmatrix} - 
u_{P_1}^{k+r}
\rnorm \\
& + \lnorm
\begin{pmatrix}
u_{P_3}^{k+r}\\
0
\end{pmatrix} - 
\begin{pmatrix}
\widehat{u}_{P_3}^{k+r} \nonumber  \\
0
\end{pmatrix} \rnorm \nonumber \\ 
&\le e_1 + \epsilon + e_2 = C_4 h + C_5 \epsilon. \label{e4}
\end{align} 

\textbf{Step 2:} Since $n$ is large enough to resolve the solution, the difference between $u_{P_2}^{k+r}$ (prolonged by padding with zeros) and $u_{P_1}^{k+r}$ should be small:
\begin{align}
e_4 = \lnorm
\begin{pmatrix}
u_{P_2}^{k+r}\\
0
\end{pmatrix} - u_{P_1}^{k+r}\rnorm \le \epsilon. \label{e5}
\end{align}

\textbf{Step 3:} Inequalities \cref{e4} and \cref{e5} give
\begin{align}
||u_{P_2}^{k+r} - u_{P_3}^{k+r}|| & \le \lnorm
\begin{pmatrix}
u_{P_2}^{k+r}\\
0
\end{pmatrix} - u_{P_1}^{k+r} \rnorm +
\lnorm  u_{P_1}^{k+r} - 
\begin{pmatrix}
u_{P_3}^{k+r}\\
0
\end{pmatrix}
\rnorm  \nonumber \\
&= e_4 + e_3 = C_4 h + C_6 \epsilon. \label{e6}
\end{align}

\begin{figure}[tbhp]
    \centering
    \includegraphics[scale=0.75]{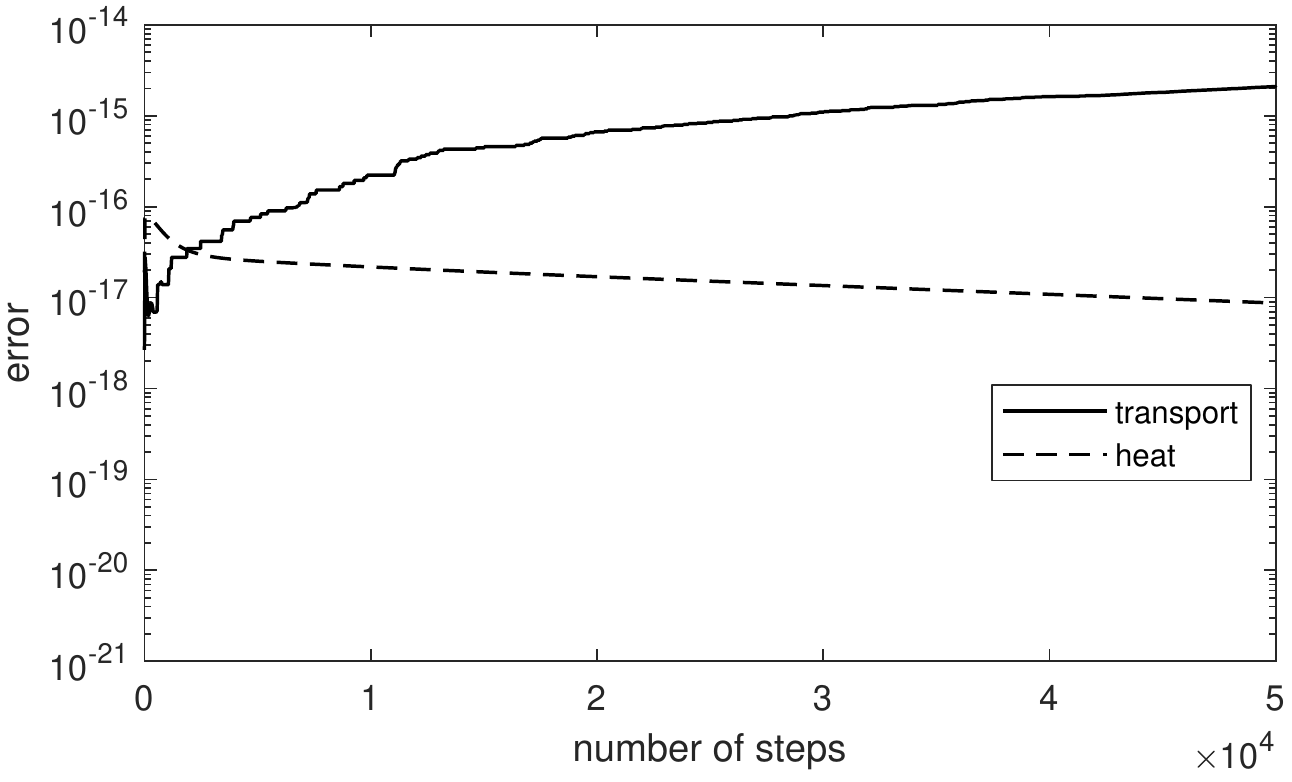}
    \caption{Difference between solutions obtained via the two approaches.}
    \label{fig:A1vsA2}
\end{figure}

The message conveyed by \cref{e6} is that the solutions obtained using Approaches 1 and 2 differ only by a quantity of $\mO(h)$ plus a multiple of $\epsilon$. In fact, the actual differences observed in all of our experiments are rather minuscule. In \cref{fig:A1vsA2}, the differences between the computed solutions via the two approaches are displayed versus the number of time steps for the one-dimensional transport equation \cref{wave1} and the heat equation \cref{heat}. To have the initial conditions and the solutions at the subsequent time steps fully resolved, we let $n = 300$ for both problems. For the one-dimensional transport equation, $4$th-order Adam-Bashforth method is used with $h = 0.1/n^2$, while for the heat equation, $3$rd-order Runge-Kutta method is used with $h = 1/n^4$. These $h$'s are chosen to stabilize the time steppers and the derivation of these restrictions is discussed in the next section. For both problems, the computed solution via the two approaches differ only by an amount of virtually machine epsilon after the first $50,000$ steps.

\section{Stability} \label{sec:stability}
\begin{figure}[tbhp]
\centering
\includegraphics[scale=0.75]{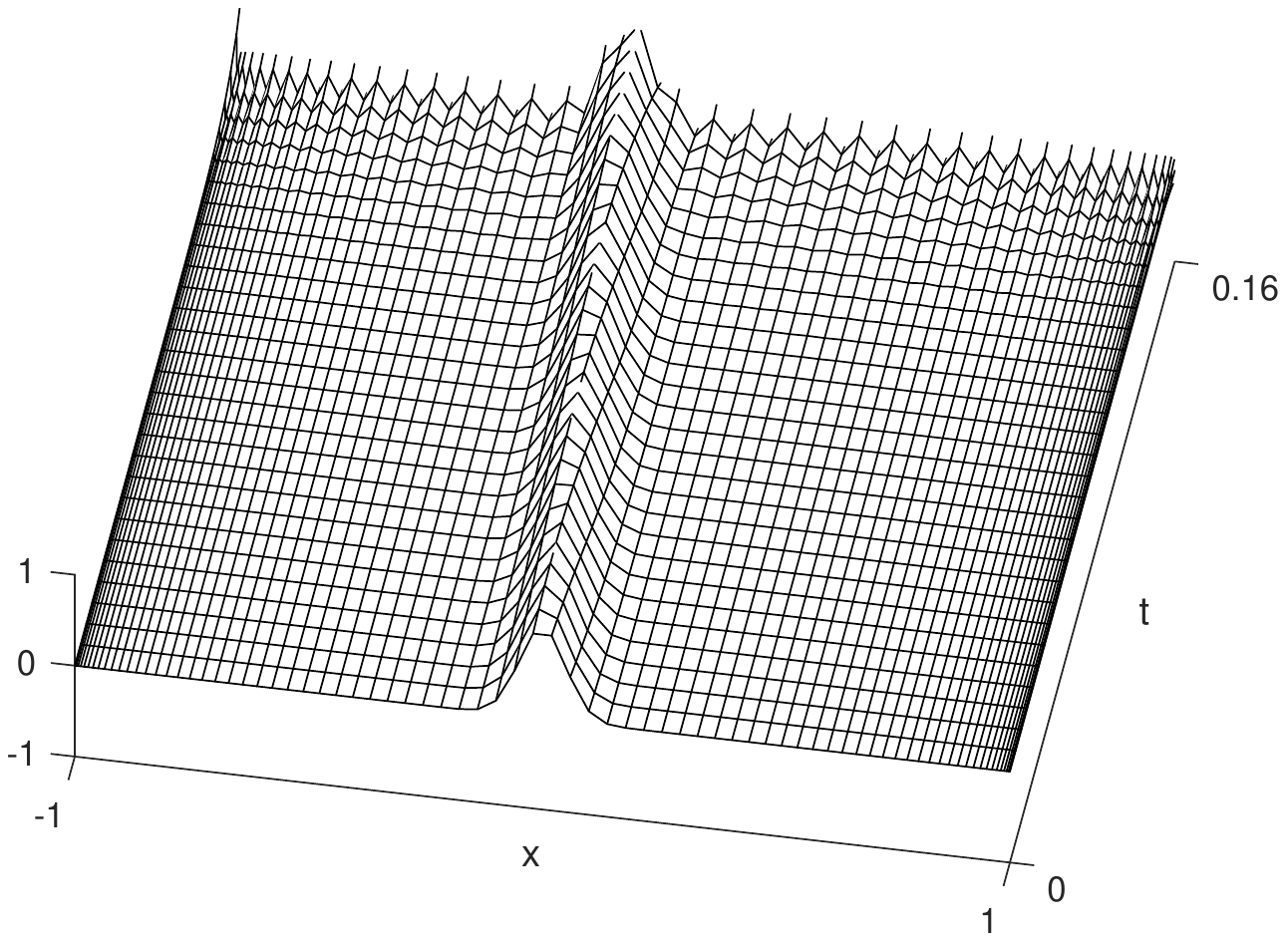}
\caption{Modal instability incited when solving \cref{wave1} with $n = 80$ and 
$h = 3.45/n^2$.}
\label{fig:unstable}
\end{figure}

The primary task now is to understand when the approaches proposed in the last 
section lead to stable calculation if a standard time stepping scheme is 
employed. In a general application of the method of lines, we consider the 
semi-discretized problem $\mT u = Au$, where $A$ is a matrix that approximates 
the spatial operator. The rule of thumb for stability is that the MOL is stable 
if the eigenvalues of $A$, scaled by the step size $h$, lie in the stability 
region of the time stepper \cite{tre1}. The same conclusion is drawn from our 
extensive experiments with the ultraspherical spectral method --- for many 
problems in the form of \cref{tdpl}, the two proposed approaches lead to 
discretization whose stability is mainly determined by the spectra of $A$. For 
example, if we apply the forward Euler method to the one-dimensional transport equation 
\cref{wave1}, both of the approaches begin to yield unstable results when $h > 
3.41/n^2$, as shown in \cref{fig:unstable}. This instability is \textit{modal} 
\cite[section 31]{tre4}, as it sets in globally and never ceases to grow --- if 
we carry on to $t = 0.3$, the unstable solution would be $\mO(10^{5})$.

Modal instability also occurs with $h > 7.2/n^4$ when the proposed approaches and the forward Euler method are used to solve the heat equation \cref{heat}. 

We now show that these restrictions are indeed attributable to eigenvalues. The discussion below will mainly be based on Approach 2 as it offers an easier form for analysis. To facilitate our discussion, we adopt the following notations in this section for the truncated version of the solution vector and the relevant operators:
\begin{align}
\begin{split}
&u = \mP_{n} \bu,\ S_{\lambda} = \mP_{n}\mS_{\lambda}\mP_{n}^{\top},\ D_{\lambda} = \mP_{n}\mD_{\lambda}\mP_{n}^{\top},\ L = \mP_{n}\mL \mP_{n}^{\top}, \\ 
&\Theta_{\lambda} = \left(\mP_{n}{\mS}_{N-1}\mP_{n+2}^{\top} \right) \left( \prod_{i=2}^{N - \lambda} \mP_{n + 2(i-1)}{\mS}_{N-i}\mP_{n + 2i}^{\top}\right),\\
&M_{\lambda} = \mP_{n+2(N-\lambda)}\mM_{\lambda}[a^{\lambda}]\mP_{n}^{\top}\text{ for }\lambda = 0, 1, \ldots, N.
\end{split} \label{PopP}
\end{align}

Let us first look at the one-dimensional transport equation \cref{wave1}. When Approach 2 is applied to \cref{wave1}, the system we end up solving can be written as
\begin{align}
S_0\mT W u=D_1 u,  \label{stwudu}
\end{align}
where, other than the temporal differential operator $\mT$, all the operators 
and the solution vector are replaced by their discretized and truncated 
counterparts. Here, $W = 
\begin{pmatrix}\vspace{0.5ex}
\begin{array}{c|c} I_{n-1} & 0 \\
\hline 
\end{array}\\
\mB\mP_{n}^{\top}
\end{pmatrix}
$ is an $n \times n$ matrix. The Dirichlet boundary condition is represented by 
the following $1 \times \infty$ functional $\mB = \begin{pmatrix}1&1&1&1& 
\cdots& \end{pmatrix}$. 

Note that the last row of \cref{stwudu} simplifies to (see \eqref{S0})
\begin{align}
\frac{1}{2} \mT \mB  \mP_{n}^{\top} u = 0. \label{wave1_bc1}
\end{align}
Since any of the multistep and Runge-Kutta methods represents $u^{k+r}$ as a linear combination of $u^{k+j}$ for $j = 0, 1, \ldots, r-1$ and the boundary condition is satisfied by $u^{k+j}$ for $j = 0, 1, \ldots, r-1$, \cref{wave1_bc1} amounts to the statement that the boundary condition is also satisfied, that is,
\begin{align}
\mB  \mP_{n}^{\top} u = 0. \label{wave1_bc2}
\end{align}

On the other hand, the top $n-1$ rows of \cref{stwudu} are
\begin{align}
S_0\texttt{(1:n-1,1:n-1)} \mT u\texttt{(1:n-1)} + S_0\texttt{(1:n-1,n)} \mT 
\mB\mP_{n}^{\top} u = D_1\texttt{(1:n-1,:)} u. \label{tulu_truncated_wave}
\end{align}
Because of \cref{wave1_bc2}, the second term on the left-hand side of \cref{tulu_truncated_wave} can be dropped and \cref{tulu_truncated_wave} coincides with the first $n-1$ rows of \cref{tulu_truncated2} for $N=1$ and $\mL = \mD$. Hence, \cref{stwudu} represents the semi-discretized system for which the largest eigenvalue(s) of $W^{-1}{S_0}^{-1}D_1$ may determine the step size of a time-stepping method for stability. The following theorem gives an upper bound for the spectral radius of $W^{-1}{S_0}^{-1}D_1$.

\begin{theorem} \label{thm:wave1}
    The spectral radius of $Q = W^{-1}{S_0}^{-1}D_1$ satisfies
    \begin{align*}
    \rho(Q) \leq (n-1)^2\sqrt{\frac{1}{3} + \frac{2}{3(n-1)^2}}.
    \end{align*}
\end{theorem}
    
\begin{proof}
    We assume that $n$ is an odd number; the proof for the even case follows analogously. Note that 
    \begin{align*}
    S_{0} = (I - B)A,
    \end{align*}
    where $A = \displaystyle \diag \left(1 , \frac{1}{2}, \frac{1}{2}, \cdots, \frac{1}{2}\right)$ and $B = \begin{pmatrix} 0_{(n-2)\times 2}&I_{n-2}\\
    0_{2 \times 2} &0_{2\times(n-2)} \end{pmatrix}$. Since $B$ is a double-shift matrix, the inverse of $S_{0}$ can be represented as a finite series 
    \begin{align*}
    S_{0}^{-1} = A^{-1}(I - B)^{-1}= A^{-1} \sum_{j = 0}^{n/2}B^{j},
    \end{align*}
    which, when spelled out, reads
    {\footnotesize{
    \begin{align} 
    S_{0}^{-1} =
    \begin{pmatrix}
      & 1\;\;&\quad&1& &1 & \,\cdots\,& 1& \\
      & &2 &\quad&2 & &\ddots&\quad&\\ 
      & & & 2& & 2&\quad& \vdots&\\
      & & & &\ddots & &\ddots &&  \\
      & & & & & \ddots&\quad&2 & \\
      & & & & &\quad&2&\quad &\\
      & & & & &\quad&\quad&\;2 & \\
    \end{pmatrix}. \label{S0Inv}
    \end{align}
    }}
    and it is easy to show
    \begin{align*}
    W^{-1} =
    \begin{pmatrix}\vspace{0.5ex}
    \begin{array}{c|c} I_{n-1} & 0 \\
    \hline 
    \end{array}\\
    -\mathcal{B}\mP_{n}^{\top}
    \end{pmatrix}.
    \end{align*}
    A simple calculation gives 
    {\footnotesize{
    \begin{align*}
    Q = \begin{pmatrix}
     &1& &3& &5& & &n-2& \\
     & &4& &8& &\cdots&2(n-3)& &2(n-1)\\
     & & &6& &10& & &2(n-2)& \\ \vspace{0.5ex}
     & & & &\ddots& &\ddots& & &\vdots\\ \vspace{0.5ex}
     & & & & &\ddots& &\ddots& &\vdots\\ \vspace{0.5ex}
     & & & & & &\ddots& &\ddots&\vdots\\ \vspace{0.5ex}
     & & & & & & &\ddots& &\vdots\\ \vspace{0.5ex}
     & & & & & & & &\ddots&\vdots\\ \vspace{0.5ex}
     & & & & & & & & &2(n-1)\\
    0&-1^2&-2^2&-3^2&-4^2\quad&-5^2\quad&\quad\cdots\quad&-(n-3)^2&-(n-2)^2&-(n-1)^2\\
    \end{pmatrix}.
    \end{align*}}}
    
    The characteristic polynomial $\det (\lambda I - Q)={\lambda}^{n} + a_{n-1}{\lambda}^{n-1}+\ldots + a_1 \lambda + a_0$ has concise expressions for the coefficients of the leading terms\footnote{For an $n \times n$ matrix $A$, $E_k(A)$ denotes the sum of $A$'s principal minor of size $k$ \cite[section 1.2]{hor}, and we use the notation $A[\alpha] = A[\alpha, \alpha]$, where $\alpha$ is a set of indices, to denote a principal submatrix of $A$ \cite[section 0.7.1]{hor}.}
    \begin{subequations}
    \begin{align*}
    a_{n-1} &= -tr(Q) = (n-1)^2,\\
    a_{n-2} &= E_2(Q) = \sum_{1\leq i\ne j \leq n}\det(Q[\{i, j\}]) = \frac{(n-1)^2}{3} [(n-1)^2-1],
    \end{align*}
    \end{subequations}
    which, by Vi\`{e}ta's theorem \cite[section 5.7]{van}, imply
    \begin{subequations}
    \begin{align*}
    \sum_{k=1}^{n}{\lambda_k} &= -(n-1)^2, \\
    \sum_{i<j}{\lambda_i}{\lambda_j} &= \frac{(n-1)^2}{3}[(n-1)^2-1],
    \end{align*}
    \end{subequations}
    where $\{\lambda_k\}_{k = 1}^n$ are the $n$ roots of $\det (\lambda I - Q)$, i.e., the eigenvalues of $Q$. We then have
    \begin{align*}
    \sum_{k=1}^{n}\lambda_k^2 = \left(\sum_{k=1}^{n}{\lambda_k}\right)^2-2\sum_{i\ne j,i<j}{\lambda_i}{\lambda_j} = \frac{(n-1)^4}{3} + \frac{2(n-1)^2}{3},
    \end{align*}
    which gives 
    \begin{equation*}
    \lvert \lambda_{max}\lvert \le \sqrt{\frac{(n-1)^4}{3} + \frac{2(n-1)^2}{3}} 
    = (n-1)^2\sqrt{\frac{1}{3} + \frac{2}{3(n-1)^2}}. 
    \end{equation*} \qed
\end{proof}

The necessary condition of the step size can be readily derived from \cref{thm:wave1}. For example, for the forward Euler method to be stable, it is required that $\lvert h \lambda_{max} + 1 \lvert \le 1$, that is,
\begin{align}
h \le \frac {3.41}{(n-1)^2}. \label{threshold}
\end{align}

\begin{figure}[tbhp]
    \centering
    \subfloat[$\veps=10^{-3},10^{-4},\ldots,10^{-17}$ in the master plot, whereas $\veps=10^{-2.75},10^{-3},\ldots,10^{-4}$ in the inset.]{\label{fig:pseigSWE(new)}\includegraphics[scale=0.51]{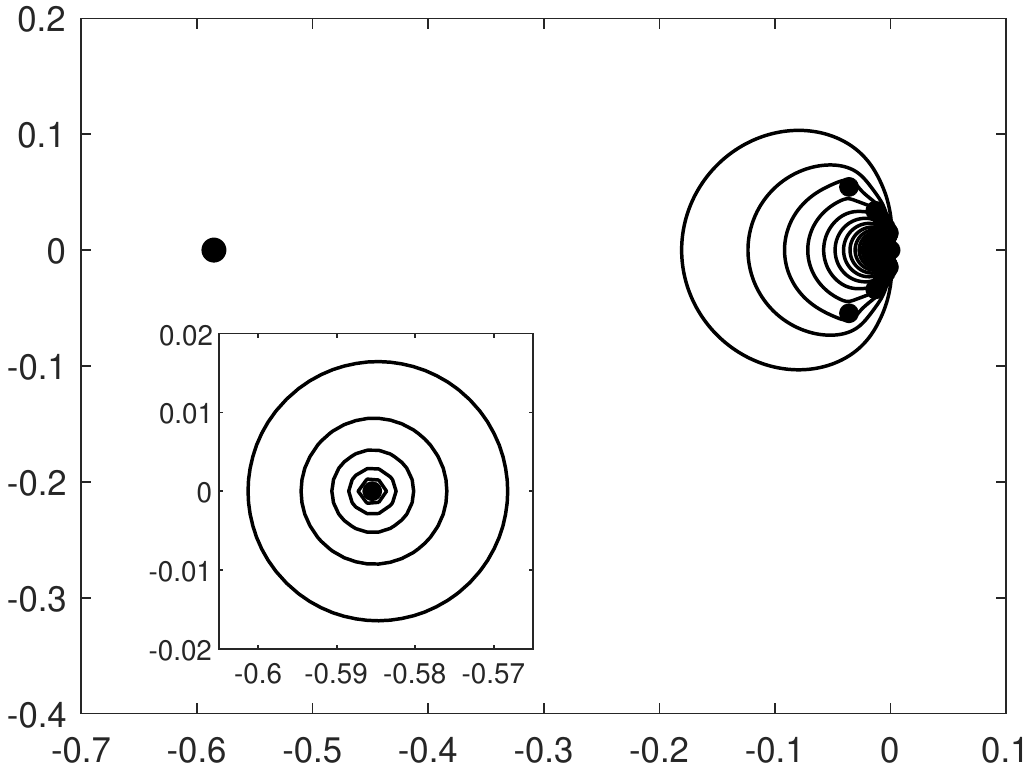}} \hspace{6mm}
    \subfloat[$\veps=10^{-3},10^{-4},\ldots,10^{-17}$.]{\label{fig:pseigSWE(old)}\includegraphics[scale=0.4234]{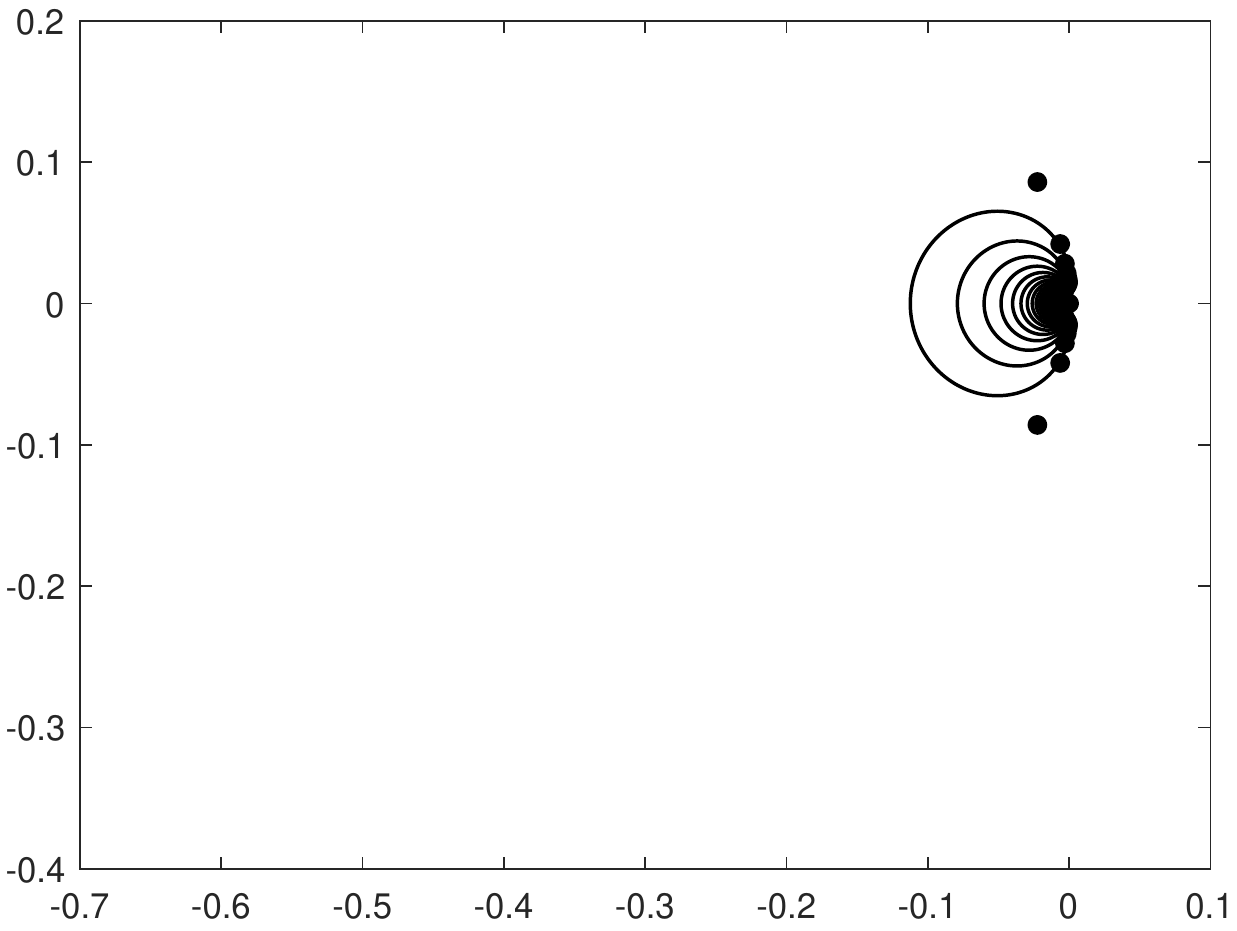}}\\
    
    \subfloat[$\veps = 10^{-2.75},10^{-3.00},\ldots,10^{-4.25}$ in both plots.] {\label{fig:pseigHE(new)}\includegraphics[scale=0.51]{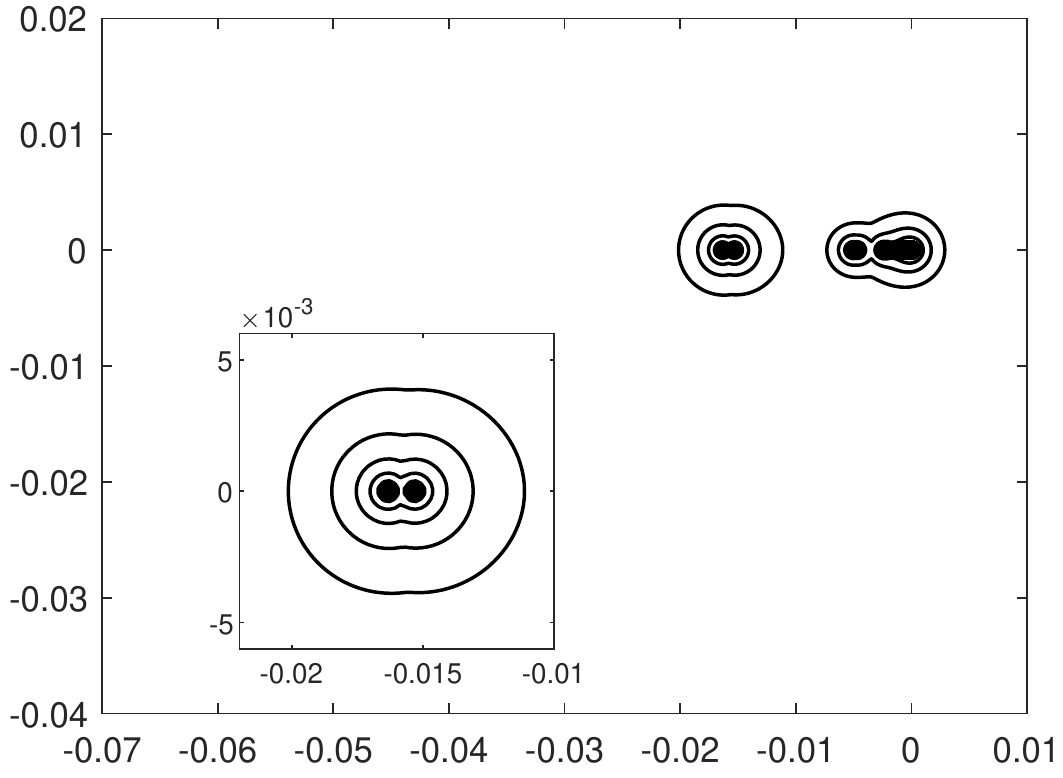}} \hspace{5mm}
    \subfloat[$\veps=10^{-2.75},10^{-3.00},\ldots,10^{-4.25}$ in both plots.]{\label{fig:pseigHE(old)}\includegraphics[scale=0.51]{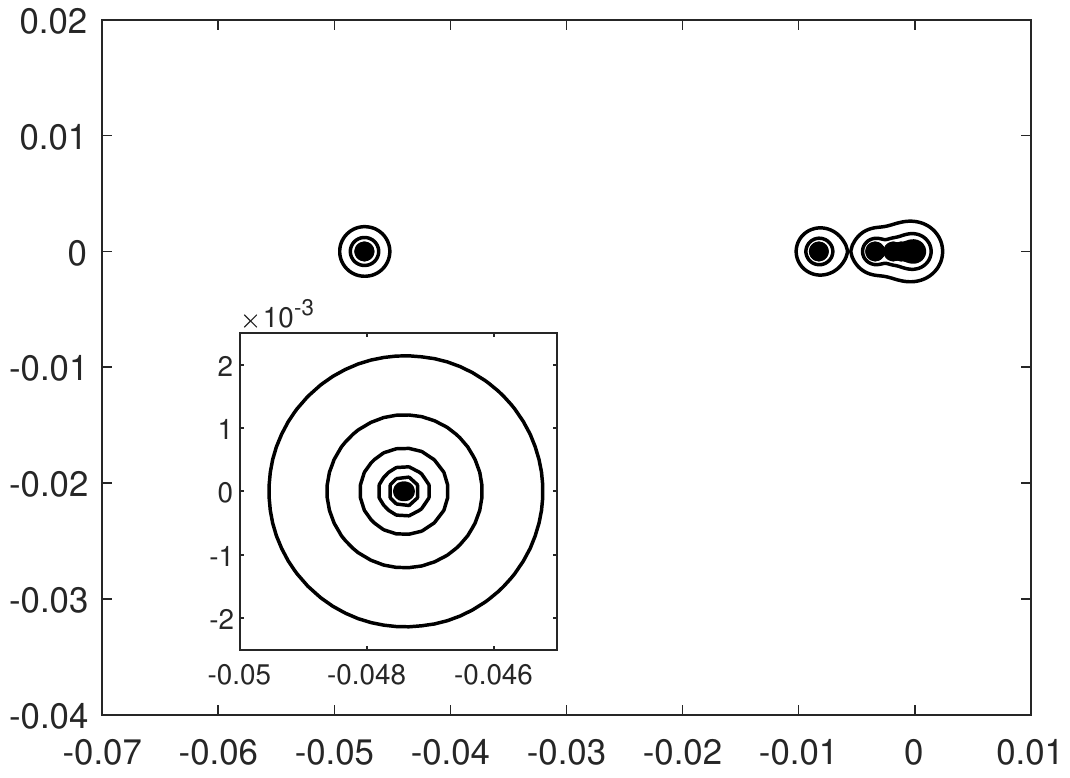}}\\
    
    \caption{The spectra and the $\veps$-pseudospectra of the spatial discretization matrices due to the ultraspherical spectral method (left panes) and the Chebyshev pseudospectral method (right panes), rescaled by $n^{-2}$ and $n^{-4}$ respectively for the one-dimensional transport equation (top panes) and the heat equation (bottom panes). The insets show close-ups in the neighborhood of the outlier(s).} \label{fig:ps}
\end{figure}

%

This is exactly the threshold beyond which we see the modal instability as in \cref{fig:unstable}. In \cref{fig:pseigSWE(new)}, the eigenvalues of $n^{-2}Q$ for $n = 64$ are plotted using dots, where the rescaling factor of $n^{-2}$ helps remove the dependence of the entries of $Q$ on the dimension $n$. The largest eigenvalue is an outlier located on the real axis that breaks away from the main cohort formed by the rest of the spectra and juts out into the left half plane. Along with the spectra, we also display the $\varepsilon$-pseudospectra in \cref{fig:pseigSWE(new)}. The pseudospectra clearly show the importance of the outlier. Although the spatial discretization matrix of the one-dimensional transport equation is nonnormal, as indicated by the pseudospectra around the main cohort, it is not far from normal in that its behavior is largely determined by the outlier. More precisely, it is only the magnitude of this outlier, not the pseudospectra around it, that matters. This can be seen from the facts that (1) the outlier is much larger in modulus than the $\varepsilon$-pseudospectra around the main cohort, (2) even in a plot for $\varepsilon$ as large as $10^{-3}$ we do not see the pseudospectra contours around the outlier, and (3) the pseudospectra around the outlier (see the close-up) consist of a few concentric circles whose radii shrink proportionally with the order of $\varepsilon$. In fact, this outlier is almost a \textit{normal eigenvalue} \cite[\S 52]{tre4} --- its condition number $\kappa(\lambda_{outlier}) \approx 9.2$ (calculated in $\infty$-norm). In contrast, the two most outlying eigenvalues in the main cohort (symmetric about the real axis) have a condition number approximately $1378.8$, which is the smallest among all the eigenvalues in the main cohort\footnote{The condition number of other eigenvalues in the main cohort could be much larger. The closer they are to the origin, the greater the condition numbers become. The eigenvalues near the origin can hardly be numerically calculated to any satisfactory accuracy due to the extremely poor conditioning.}. That is, the eigenvalues in the main cohort are nonnormal or significantly so. In a word, the outlier is of physical significance, and it is this outlier that restricts the step size of a time stepper with a bounded stability region. Hence, we can say that for the proposed approaches the stability of a time marching scheme, when applied to \cref{wave1}, is mainly determined by the spectra.

For comparison purposes, we show in \cref{fig:pseigSWE(old)} the spectra and 
pseudospectra of the rescaled first-order differentiation matrix from the 
Chebyshev pseudospectral method. As the largest eigenvalues in 
\cref{fig:pseigSWE(old)} are smaller than the outlier in 
\cref{fig:pseigSWE(new)}, one might think that the Chebyshev pseudospectral 
method is superior to the ultraspherical spectral method. This is, in fact, not 
the case. First, it is not true that ultraspherical spectral method always 
results in greater spectral radius than the collocation pseudospectral method, 
e.g., the second-order differentiation operator (see below), or if boundary 
conditions of other types are enforced (see \cref{thm:bc}). Second, the 
Chebyshev pseudospectral method, in this particular case, allows a step size 
only of a constant times larger than the ultraspherical spectral method, since 
for both methods the spectral radius is $\mO(n^2)$. Third, the ultraspherical 
spectral method is cheaper stepwise than the Chebyshev pseudospectral method 
(see \cref{sec:cost}), thanks to its sparsity structures.

Now, we turn to the heat equation \cref{heat} which features the spatial differentiation of a second order. Applying Approach 2 to the heat equation \cref{heat}, but leaving the temporal operator non-discretized, gives
\begin{align}
S_1 S_0 \mT H u = D_2 u, \label{ssthudu}
\end{align}
where $H = 
\begin{pmatrix}\vspace{0.5ex}
\begin{array}{c|c} I_{n-2} & 0 \\
\hline 
\end{array}\\
\mB\mP_{n}^{\top}
\end{pmatrix}$, and the functional $\mB = 
\begin{pmatrix}
1&1&1&1& \cdots \\
1&-1&1&-1& \cdots
\end{pmatrix}$ represents the Dirichlet boundary conditions in Chebyshev space. 
The equivalence of \cref{ssthudu} and the Approach 2 discretization follows 
exactly the same reasoning for that of \cref{stwudu} where $W$ is involved instead. 

The following theorem gives a bound on the spectral radius of $H^{-1}S_0^{-1}S_{1}^{-1}D_2$.

\begin{theorem} \label{thm:heat}
The spectral radius of $G = H^{-1}{S_0}^{-1}{S_1}^{-1}D_2$ is bounded by
\begin{align*}
\rho(G) \leq \frac{2}{3}n(n-2)(n-1)^2.
\end{align*}
\end{theorem}
\begin{proof}
$S_0^{-1}$ and $D_2$ are given by \cref{S0Inv} and \cref{opDiff}, respectively, 
and $S_{1}^{-1}$ can be derived analogously to $S_0^{-1}$. Also, we have 
$H^{-1} = 
\begin{pmatrix}\vspace{0.5ex}
\begin{array}{c|c} I_{n-2} & 0 \\
\hline 
\end{array}\\
H'
\end{pmatrix}$,
where $H' = \begin{pmatrix}
     &-1& &-1& \ldots&0&\frac{1}{2}&-\frac{1}{2}\\
    -1& &-1& &\ldots&-1&\frac{1}{2}&\frac{1}{2}\\
    \end{pmatrix}.$
It can be shown that $G_{n n} = \displaystyle -\frac{2}{3}n(n-2)(n-1)^2$ is the entry with the largest magnitude. Noting this, we can further show that for any $\lambda < -\displaystyle \frac{2}{3}n(n-2)(n-1)^2$ the determinant $\det (\lambda I - G) \neq 0$. Hence, all eigenvalues of $G$ are smaller than $\displaystyle \frac{2}{3}n(n-2)(n-1)^2$ in modulus. 
\qed
\end{proof}

The spectra and the pseudospectra of $n^{-4}G$ are shown in \cref{fig:pseigHE(new)}, where the eigenvalues are lined up on the real axis, due to the parity of the order of the spatial differentiation. Once again, there are (two) outliers which detach themselves from the rest of the spectra and reside far in the left half plane. Though we can see the pseudospectra for both the outliers and the rest of the spectra, the relatively large $\varepsilon$, the relatively small scale of the axes, the shape of the pseudospectra contours around the outliers, and the fact that the condition numbers of these two outliers are small (both approximately $2.3$) suggest that the outliers are physically significant, governing the behavior of the matrix\footnote{In fact, the rest of the eigenvalues are all normal with $\mO(1)$ condition numbers.}. Therefore, the outlier of largest modulus determines the maximum step size if a time stepper with bounded stability region is used. The spectra and the pseudospectra of the rescaled Chebyshev second-order differentiation matrix are shown in \cref{fig:pseigHE(old)} for comparison.

Again, we can derive from \cref{thm:heat} a threshold value below which the 
step size of the time marching scheme leads to a stable solution to 
\cref{heat}. Although this bound is not sharp, i.e., a step size bigger than 
this value may well stabilize the computation, the key point is not missed --- 
the largest eigenvalue of the ultraspherical discretization matrix behaves like 
$\mO(n^4)$. This echoes \cite{wei}, which gives a same result for the 
second-order pseudospectral differentiation matrix. Such an agreement is not a 
coincidence. Furthermore, the last two theorems suggest that the largest eigenvalues of the $N$th order spatial differentiation operator, when truncated to $n \times n$ and converted back to Chebyshev space, 
scale like $\mO(n^{2N})$, the same as in the Chebyshev pseudospectral 
methods\footnote{It is well known that the largest eigenvalues of the $N$th 
order Chebyshev pseudospectral differentiation matrix scale like $\mO(n^{2N})$. 
Surprisingly, however, this assertion is not found in the literature and it 
seems that no one has ever given a proof of it.}. Indeed, this is exactly what 
we show in \cref{thm:rho} below. To do so, we look at \cref{tulu_truncated2}. 
Inverting the product of the conversion matrices on the left-hand side of 
\cref{tulu_truncated2} gives
\begin{align}
\mT u = \left(S_{N - 1} S_{N - 2} \ldots S_{0}\right)^{-1}Lu, \label{tulu_truncated4}
\end{align}
where each conversion matrix is truncated exactly before the inversion. The following lemma gives an upper bound for the norm of $S_{\lambda}^{-1}$.  
\begin{lemma} \label{thm:S}
For $\lambda = 1, 2, \ldots $,  $\lnorm S_{\lambda}^{-1}\rnorm \leq C_{\lambda}n^2$ for some constant $C_{\lambda}$ and $\lnorm S_0^{-1}\rnorm \leq n$.
\end{lemma}
\begin{proof}

Following a derivation analogous to the one for $S^{-1}_{0}$, we find
\begin{align*}
S_{\lambda}^{-1} =
\begin{pmatrix}
1& &1 & &1 & & \dots& 1&  \\
 & \frac{\lambda + 1}{\lambda} & & \frac{\lambda + 1}{\lambda} & & \frac{\lambda + 1}{\lambda} & & &  \\
 &  & \frac{\lambda + 2}{\lambda}&  & \frac{\lambda + 2}{\lambda} & &\ddots & \vdots&  \\
 & & & \frac{\lambda + 3}{\lambda} & & \frac{\lambda + 3}{\lambda}& &\vdots & \\
 & & & &\ddots & &\ddots & & \\
 & & & & &\ddots & &\frac{\lambda + n-3}{\lambda}  \\
 & & & & & &\frac{\lambda + n-2}{\lambda} & \\
 & & & & & & &\frac{\lambda + n-1}{\lambda} & \\
\end{pmatrix}. 
\end{align*}

Hence, we have
\begin{align*}
\lnorm S_{\lambda}^{-1}\rnorm &= \max_i \left(\frac{n + 1}{2}, \frac{n - 1}{2} \frac{\lambda + 1}{\lambda} , \cdots ,\left(\frac{n +1}{2} - \left \lceil \frac{i}{2} \right \rceil \right) \frac{\lambda + i}{\lambda}, \cdots  \frac{\lambda + n - 1}{\lambda}\right) \\
&\leq C_{\lambda}n^2,
\end{align*}
and $\lnorm S_0^{-1}\rnorm \leq n$ follows from \cref{S0Inv}.\qed
\end{proof}

Now we are in a position to bound the norm of the matrix on the right-hand side of \cref{tulu_truncated4}. 
\begin{lemma} \label{thm:SPLP}
Suppose that each of $M_{\lambda}$ is of a finite bandwidth independent of the degrees of freedom $n$ for $\lambda = 0, 1, \ldots, N$, then
\begin{align}
\lnorm \begin{pmatrix}
S_{N-1} \dots S_0
\end{pmatrix}^{-1}
L\rnorm \le Cn^{2N} \label{SL}
\end{align}
for some constant $C$.
\end{lemma}

\begin{proof}
From \cref{opDiff}, it is easy to see that  
\begin{align*}
\lnorm D_{\lambda}\rnorm  \leq C_{N}n
\end{align*}
for all $\lambda$. 

Since $M_{\lambda}[a^{\lambda}]$ has a finite bandwidth, $\lnorm M_{\lambda}[a^{\lambda}]\rnorm $ is bounded by a constant regardless the dimension $n$. Similarly, this is the case for $\lnorm S_{\lambda}\rnorm $ for all $\lambda$. 

By the triangle inequality and the submultiplicativity of matrix norms, it follows from \cref{PopP} that
\begin{align*}
\lnorm L\rnorm \le n \left(\lnorm M_N\rnorm +\sum_{\lambda=1}^{N-1} \lnorm \Theta_{\lambda} M_{\lambda}\rnorm \right) \le C_L n,
\end{align*}
for some $C_L$ and this, along with \cref{thm:S}, gives \cref{SL}.
\qed
\end{proof}

A direct consequence of \cref{thm:SPLP} is an upper bound for the spectral radius of the matrix on the right-hand side of \cref{tulu_truncated4}.
\begin{theorem} \label{thm:rho}
When Approach 2 is used for solving \cref{tdpl} where $\mL$ is an $N$th order 
differential operator with smooth variable coefficients given by \cref{opL}, 
there exists a constant $C$ independent of the degrees of freedom $n$ for the 
spatial discretization such that
\begin{align}
\rho(S_0^{-1} S_1^{-1} \ldots  S_{N-1}^{-1} L) \leq Cn^{2N}. \label{srbound}
\end{align}
\end{theorem}
\begin{proof}
The smoothness of the variable coefficients implies finite bandwidth for 
each $M_{\lambda}$. Hence, this is a standard result led to by \cref{SL} which 
can be found, for example, in \cite[Theorem 5.6.9]{hor}.
\qed
\end{proof}

\begin{figure}[tbhp]
  \centering
  \subfloat[$N=3$]{\label{fig:ctuxxx}\includegraphics[scale=0.41]{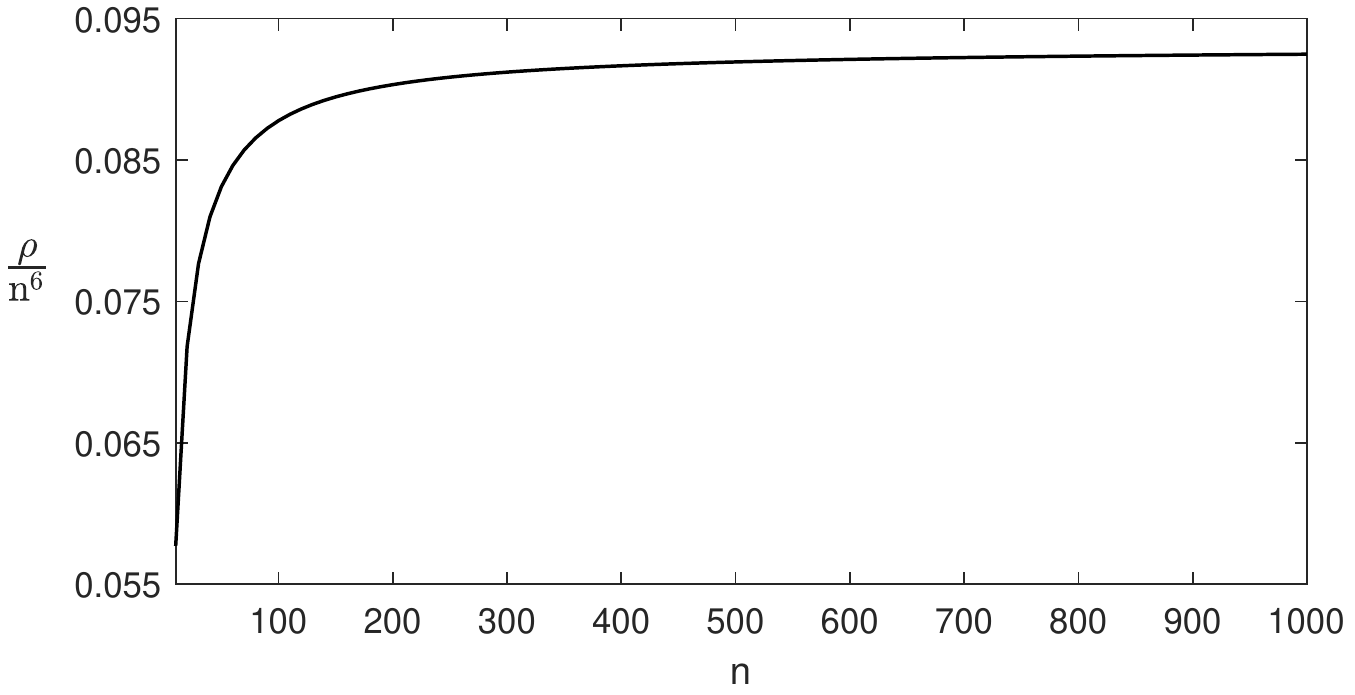}} \hspace{4mm}
  \subfloat[$N=4$]{\label{fig:ctuxxxx}\includegraphics[scale=0.41]{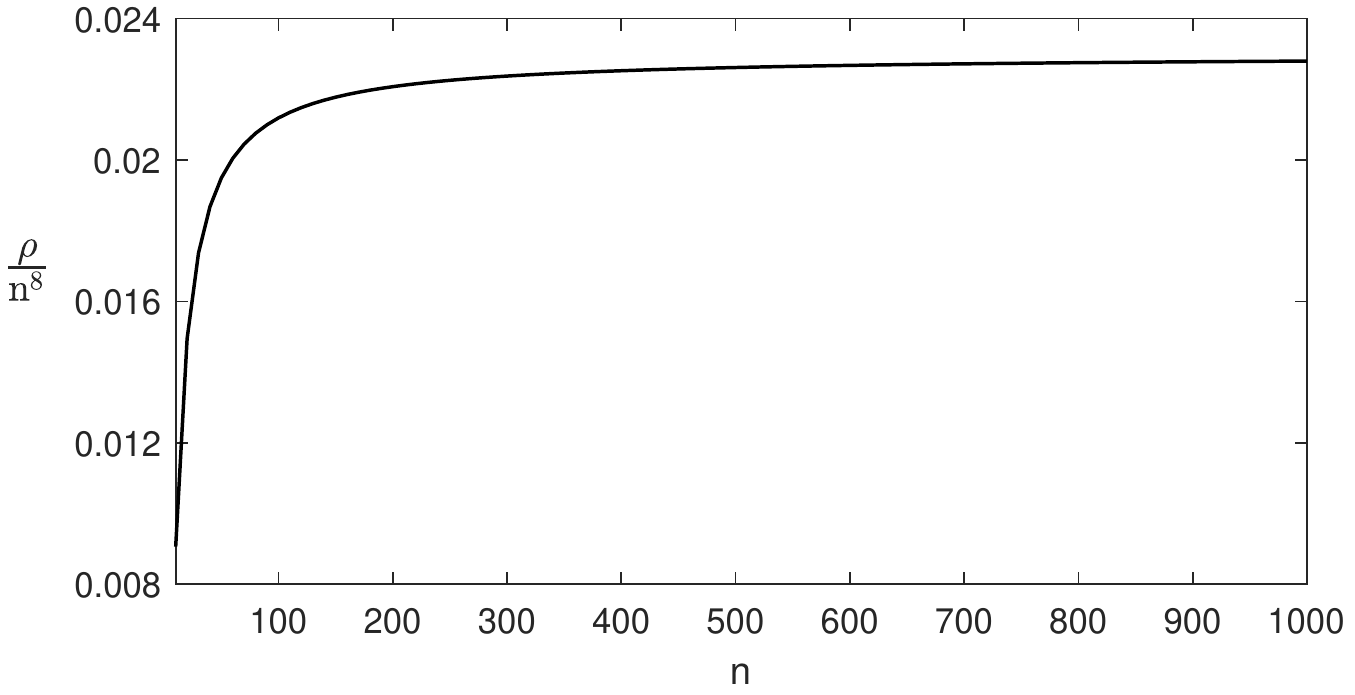}} 
  \caption{The spectral radius, normalized by $n^{-2N}$, of the $n \times n$ spatial discretization matrices for $N$th-order differentiation operators versus $n$.}
  \label{fig:cthighorder}
\end{figure}

\cref{thm:rho} is numerically verified for the cases of $N = 3, 4$ in \cref{fig:cthighorder}, where the spectral radius of $N$th-order differentiation matrices is normalized by $n^{-2N}$ and plotted versus different $n$. It can be seen that the normalized spectral radii indeed tend to be a constant.

In fact, the $\varepsilon$-pseudospectra radius of the matrix on the right-hand side of \cref{tulu_truncated4} is bounded by the same quantity plus $\varepsilon$.
\begin{theorem}
If the assumption holds as in \cref{thm:rho}, 
\begin{align}
\rho_{\eps}(S_0^{-1} S_1^{-1} \ldots S_{N-1}^{-1} L) \leq Cn^{2N} + \eps, \label{psbound}
\end{align}
where the constant $C$ is the one given in \cref{srbound}.
\end{theorem}
\begin{proof}
For any $\lnorm E \rnorm \leq \eps$,
\begin{align*}
\rho(S_0^{-1} S_1^{-1} \ldots  S_{N-1}^{-1} L + E) \leq \lnorm S_0^{-1} S_1^{-1} \ldots  S_{N-1}^{-1} L \rnorm + \lnorm E \rnorm \leq Cn^{2N} + \eps,
\end{align*}
which, by the second definition of pseudospectra \cite[\S 2]{tre1}, gives \cref{psbound}.
\qed
\end{proof}

The bounds in the last two theorems give the worst case scenario of how the spectra and the pseudospectra scale with $n$ for $N$. If we use the quantity $n^{2N}$ as a guidance for choosing the step size, the stability is guaranteed. 

Since the largest eigenvalue(s) of a spatial discretization matrix also grows like $\mO(n^{2N})$ for the Chebyshev pseudospectral method, the ultraspherical spectral method and the Chebyshev pseudospectral method roughly tie in terms of the largest step that can be taken for a time marching scheme with a bounded stability region. The fact that the largest eigenvalues match for these two methods can also be seen by premultiplying both sides of \cref{tulu_truncated4} by an inverse discrete cosine transform (iDCT) matrix and ignoring the first and last rows, as this reproduces the discretization led to by the Chebyshev pseudospectral method \cite[chapter 10]{tre1}. Because the iDCT matrix is unitary, the norms of the spatial discretization matrices due to these two methods should be roughly same.

\begin{remark} \label{thm:bc}
Our discussion in this section is based on the one-dimensional transport 
equation and the heat equation subject to homogeneous Dirichlet boundary 
conditions. However, the use of homogeneous Dirichlet boundary conditions is 
unimportant. Although other boundary conditions may lead to different constants 
in bounds such as \cref{threshold}, it would not change the main result given 
by \cref{thm:rho}. In addition, homogeneous Dirichlet boundary conditions were 
adopted in the study of the collocation-based pseudospectral methods \cite{got,tre1,tre4,wei}. It is for comparative purposes that the use of the same 
boundary conditions seems natural.
\end{remark}

\section{Error} \label{sec:error}

The error in the computed solution of PDEs comes mainly from two sources: discretization and rounding, where the former, in the present context, consists of those in space and time. That is,
\begin{align*}
\text{Total error} ~~~=~~~ 
\begin{matrix}
\text{spatial}\\
\text{discretization}\\
\text{error}
\end{matrix} ~~~+~~~ 
\begin{matrix}
\text{temporal}\\
\text{discretization}\\
\text{error}
\end{matrix}
 ~~~+~~~ 
\begin{matrix}
\text{rounding}\\
\text{error}
\end{matrix}.
\end{align*}

Like any other spectral method, the ultraspherical spectral method offers 
spectral accuracy, that is, the accuracy is limited not by the order of the 
discretization, but by the smoothness of the solution being approximated. When 
the degrees of freedom are sufficiently large, the solution can be adequately 
resolved in space, thereby bringing no spatial discretization error. The 
temporal discretization error introduced by the standard time marching schemes 
is usually of an algebraic order and its quantification and analysis can be 
found in standard texts like \cite{asc,but,hai}. When the time step is small 
enough, the temporal discretization error can essentially be restricted to or 
below the level of machine epsilon. Hence, it is possible to completely 
annihilate the discretization error and this is a common working paradigm 
adopted by spectral methods for PDEs. This way, one is only left with the 
errors introduced by rounding. We now give an analysis of the rounding error 
for the proposed method, assuming the discretization error is absent.   

We consider the iterative model
\begin{align}
AU^{k+1} = BU^{k}, \label{iterative}
\end{align}
which can be deemed as the prototype of the discretized systems obtained by the proposed method. Here, $U^k$ and $U^{k+1}$ are the computed solutions at two successive steps\footnote{For a $r$-step linear multistep method, such a relation can be derived by forming $A$ and $B$ as $rn \times rn$ block matrices and $U^k$ and $U^{k+1}$ as vectors that incorporate the computed solution at $r$ successive time steps.}.

The key to our analysis is the quantity
\begin{align*}
\Delta^{k+1} = \udU^{k+1} - A^{-1}B\udU^k,
\end{align*}
where $\udU^k$ and $\udU^{k+1}$, stored as floating point numbers, are the computed solutions at $k$th and $(k+1)$th step, respectively. Here, the matrix $A^{-1}B$ is assumed to be exact, not in its floating point representation, so that $\Delta^{k+1}$ quantifies the  amount of error introduced by rounding at a single step. We shall find an upper bound for its magnitude as follows.

\begin{align*}
\lnorm \Delta^{k+1}\rnorm  = \lnorm \udU^{k+1} - A^{-1}B\udU^k\rnorm  = 
\lnorm fl(\udA^{-1}\udB \,\udU^k )- A^{-1}B\udU^k\rnorm ,
\end{align*}
where $fl(x)$ denotes the function producing the closest floating point approximation to a given number $x$. There exists $\eps$ with $\lvert \eps \lvert \leq \eps_{mach}$ such that $fl(x) = x(1+\eps)$ \cite{ove}. Here, $\eps_{mach}$ is the machine epsilon and in IEEE double precision arithmetic $\eps_{mach}$ is $2^{-53} \approx 1.11 \times 10^{-16}$. Hence,
\begin{align*}
\lnorm \Delta^{k+1}\rnorm  &= \lnorm \udA^{-1}\udB \,\udU^k (1 + \eps) - A^{-1}B\udU^k\rnorm  \\ 
&\le \lnorm \udA^{-1}\udB - A^{-1}B\rnorm \,\lnorm \udU^k\rnorm  + \eps\lnorm \udA^{-1}\udB\rnorm \,\lnorm \udU^k\rnorm  \\ &\le 
\lnorm \udA^{-1}\rnorm \, \lnorm \udB - B \rnorm \, \lnorm \udU^k\rnorm + \lnorm  B\rnorm \, \lnorm \udA^{-1} - A^{-1} \rnorm \, \lnorm \udU^k\rnorm + \eps \lnorm \udA^{-1}\udB \rnorm \,\lnorm \udU^k\rnorm,
\end{align*}

By Theorem 2.3.9 in \cite{wat}, we have 
\begin{align*}
\lnorm \udA^{-1} - A^{-1} \rnorm \le C_1 \eps_{mach}, 
\end{align*}
where $\displaystyle C_1 = \frac{\lnorm A^{-1} \rnorm + \kappa(A)}{1 - \eps_{mach}\lnorm A^{-1} \rnorm} \lnorm A^{-1} \rnorm $ and $\kappa(A)$ is the condition number of $A$ in the infinity norm. A little algebraic work gives
\begin{align}
\lnorm \Delta^{k+1}\rnorm &\le C_2 \lnorm \udU^k\rnorm \epsilon_{mach}, \label{Delta}
\end{align}
where $\displaystyle C_2 = \lnorm A^{-1} \rnorm \left( n + \frac{\lnorm A^{-1} \rnorm + \kappa(A)}{1 - \eps_{mach}\lnorm A^{-1} \rnorm} \lnorm B \rnorm + \lnorm B \rnorm \right)$.

\begin{figure}[bhp]
    \centering
    \subfloat[ ]{\label{fig:err_wave}\includegraphics[scale=0.412]{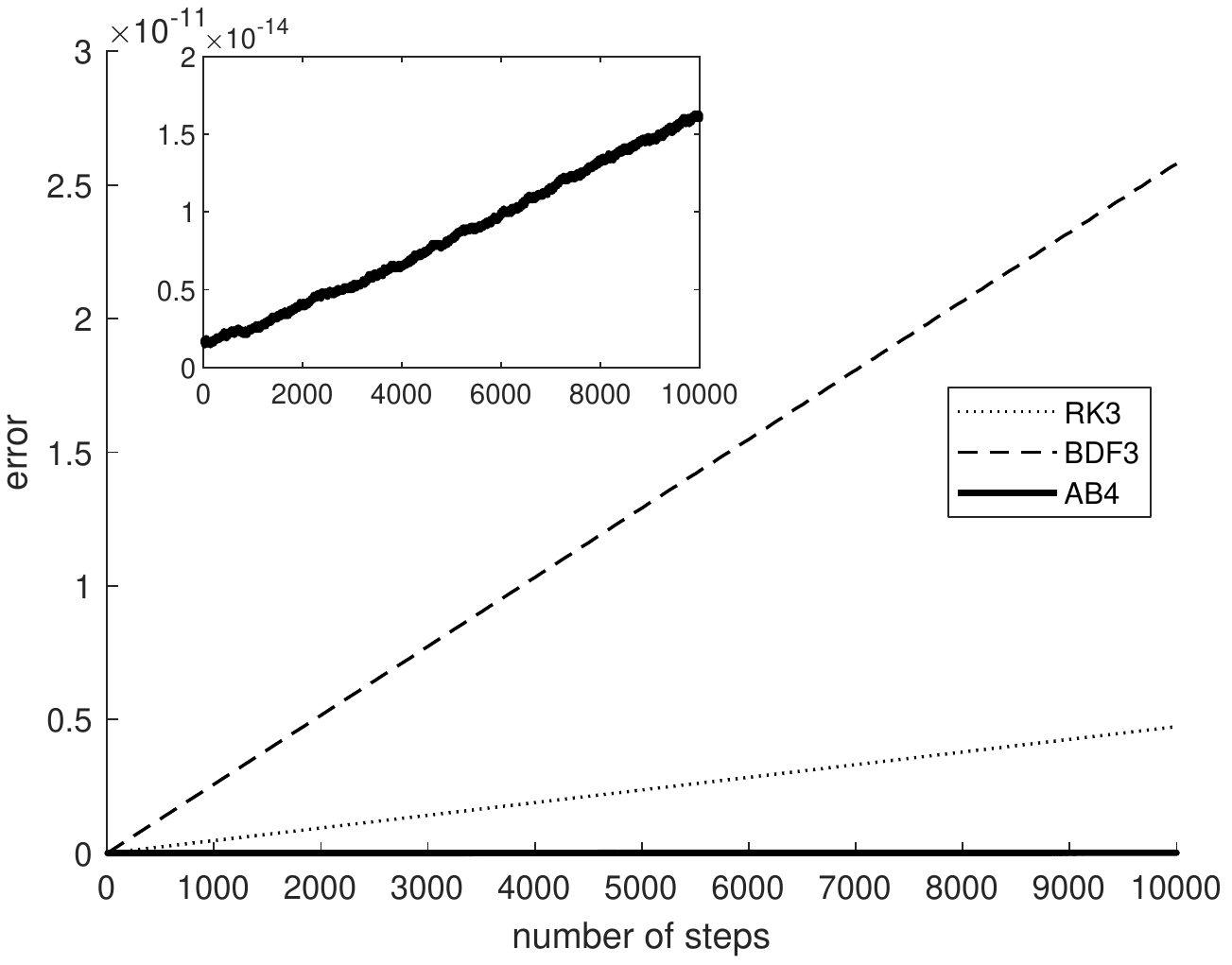}} \hspace{5mm}
    \subfloat[ ]{\label{fig:err_waveold}\includegraphics[scale=0.412]{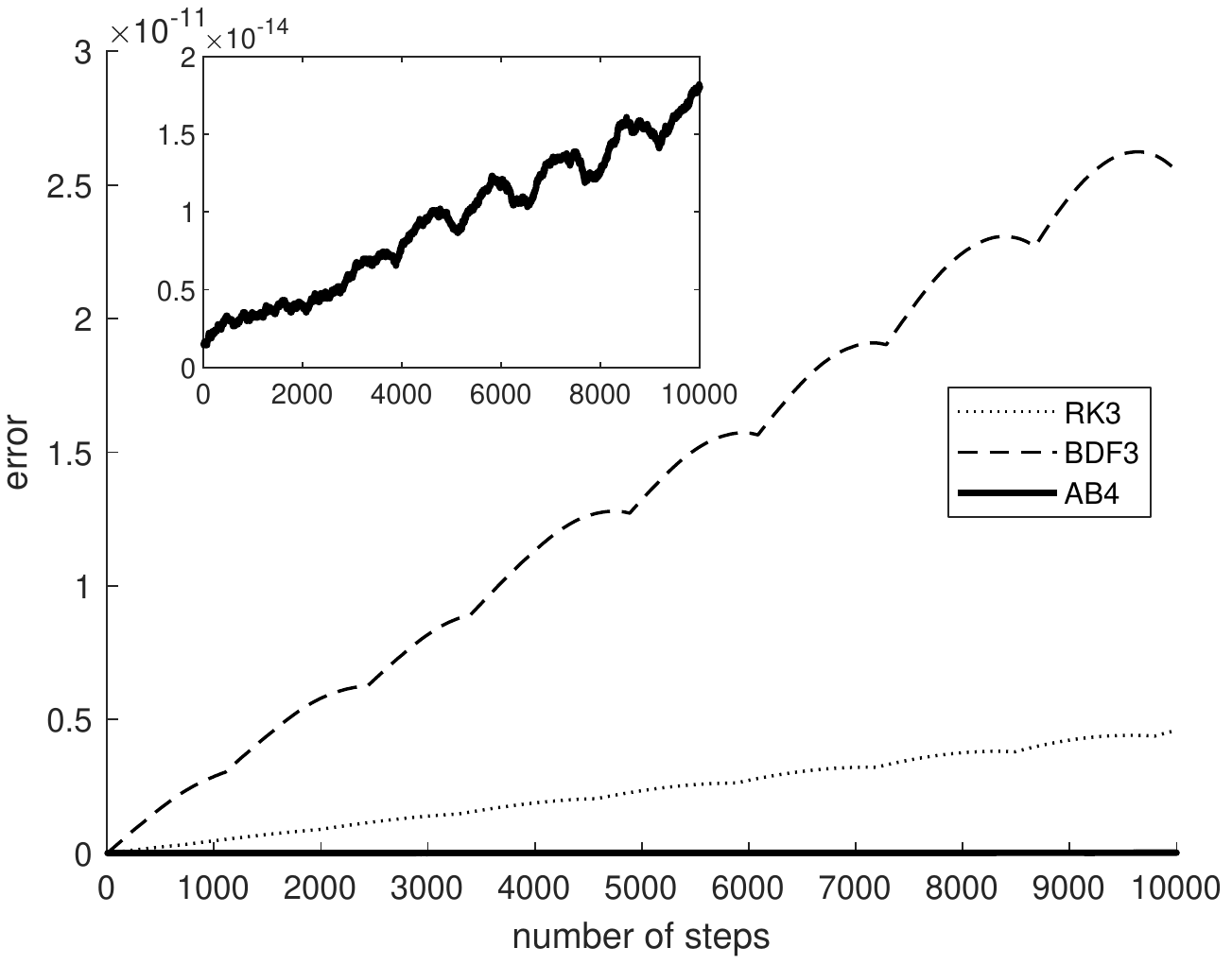}} 
    
    \subfloat[ ]{\label{fig:err_heat}\includegraphics[scale=0.43]{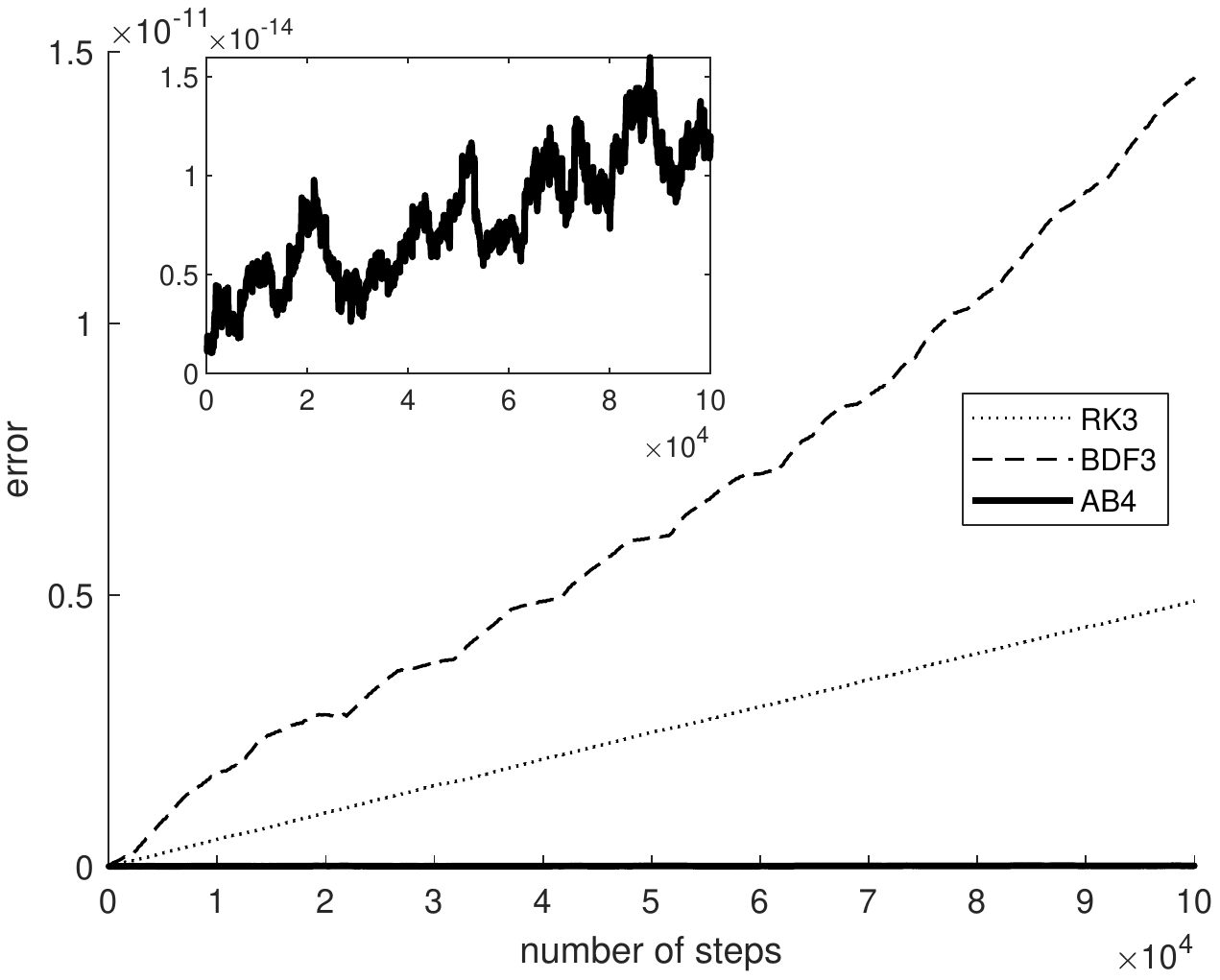}}
    \hspace{6mm}
    \subfloat[ ]{\label{fig:err_heatold}\includegraphics[scale=0.43]{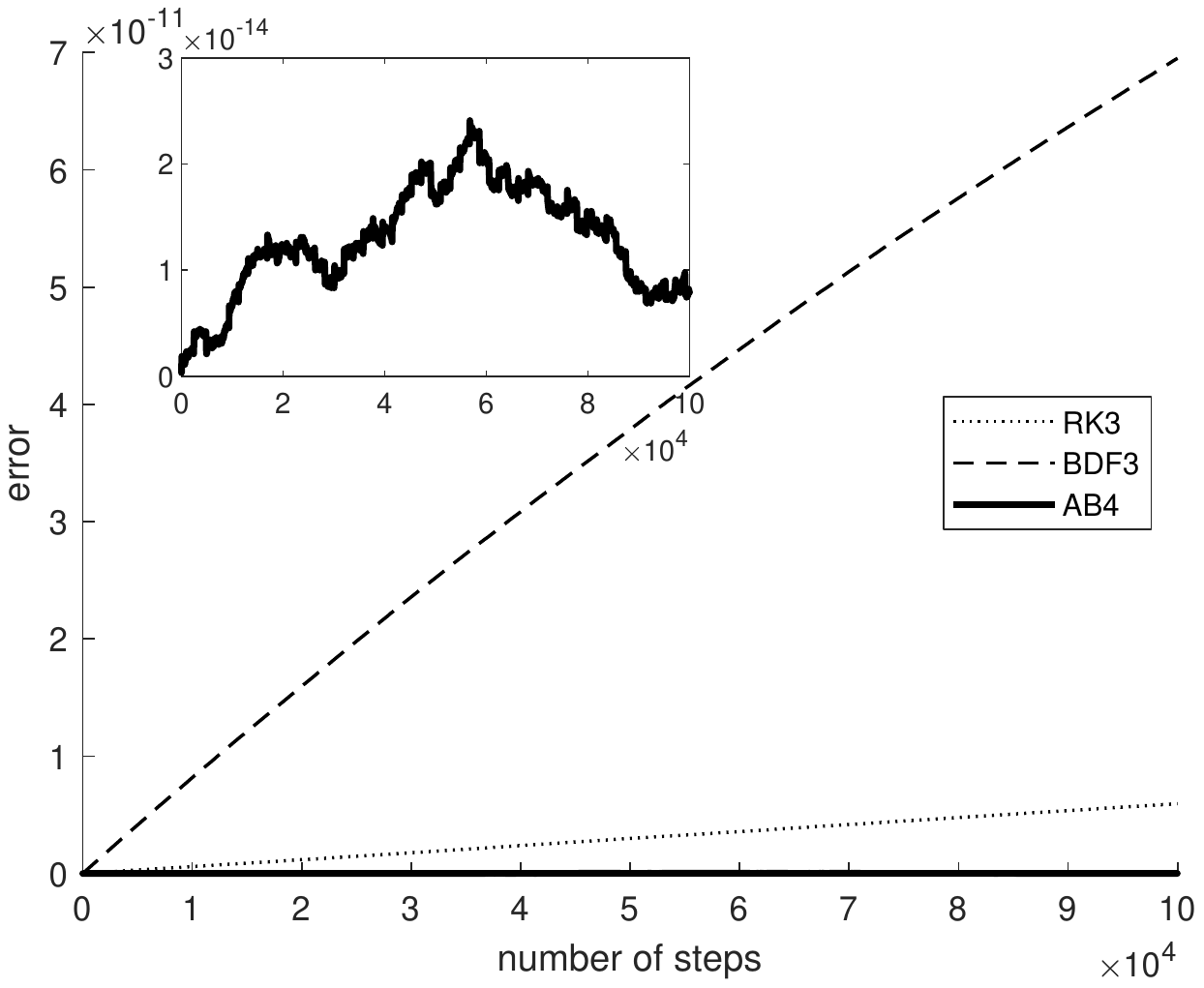}}
    \caption{The growth of the rounding error when solving \cref{wave1} (top panes) and \cref{heat} (bottom panes) using ultraspherical spectral method (left panes) and Chebyshev pseudospectral method (right panes).}
    \label{fig:error}
\end{figure}

Now we add up the error introduced in each and every step to have the accumulated error at $(K+1)$th step in the form of a discrete convolution 
\begin{align*}
E^{K+1} = \sum_{j = 1}^{K+1}\left(A^{-1}B\right)^{K+1-j}\Delta^{j},
\end{align*}
whose magnitude can be bounded as
\begin{align}
\lnorm E^{K+1}\rnorm &= \lnorm \sum_{j = 1}^{K+1}\left(A^{-1}B\right)^{K+1-j}\Delta^{j}\rnorm \\  
& \le (K+1) \sup_{0\le r\le K}\lnorm \left(A^{-1}B\right)^{r}\rnorm \sup_{1\le j\le K+1}\lnorm \Delta^{j}\rnorm \nonumber \\ 
&\le (K+1) \sup_{0 \le r \le K} \lnorm \left(A^{-1}B\right)^{r}\rnorm \sup_{1\le j\le K+1} \lnorm \underline U^j \rnorm C_2 \epsilon_{mach} \nonumber\\ 
&= C_3(K+1)\eps_{mach},
\label{accumulated}
\end{align}
where we have used \cref{Delta} to come up with the constant coefficient
{\footnotesize{
\begin{align*}
C_3 = \sup_{0\le r\le K}\lnorm \left(A^{-1}B\right)^{r}\rnorm\sup_{1\le j\le K+1}\lnorm \udU^j\rnorm \lnorm A^{-1} \rnorm \left( n + \frac{\lnorm A^{-1} \rnorm + \kappa(A)}{1 - \eps_{mach}\lnorm A^{-1} \rnorm} \lnorm B \rnorm + \lnorm B \rnorm \right),
\end{align*}}}
independent of $K$.

What \cref{accumulated} shows is that the accumulated error grows at worst \textit{linearly} with the number of the time steps. We solve the one-dimensional transport equation \cref{wave1} and the heat equation \cref{heat} using three different time marching schemes, i.e., RK3, AB4, and BDF3, and compare the computed solution with the exact solution. Sufficiently large $n$ and small enough $\Delta t$ are used so that there is no discretization error and the observed error is solely due to rounding. The error is plotted in \cref{fig:error} (left panes) to show its growth versus the number of time steps.  

As shown in \cref{fig:err_wave}, the errors grow exactly linearly for all three methods. The error of the AB4 method is relatively negligible compared to those of RK3 and BDF3, so its curve is indistinguishable from the x-axis in the same plot. The inset shows the linear growth of the error of AB4 using a different y-scale. The results shown in \cref{fig:err_heat} look similar, only except that the error curves are somewhat more oscillatory, especially in the inset plot for AB4.

The different slopes of the curves are attributed to $C_3$ in \cref{accumulated}. In fact, our calculation shows that it is the factor $\sup_{0\le r\le K}\lnorm \left(A^{-1}B\right)^{r}\rnorm$ that really makes a difference for these three methods. For example, this quantity is $3.9 \times 10^{10}$, $8.9 \times 10^6$, and $5.0$ for BDF3, RK3, and AB4, respectively. Note that this factor is partly attributed to our analysis on the norms of the spatial discretization matrices in \cref{sec:stability} and the Kreiss matrix theorem \cite[Chapter 18]{tre4}.

What we also show in \cref{fig:error} (right panes) is how the rounding errors grow when Chebyshev pseudospectral method is used to solve \cref{wave1} (\cref{fig:err_waveold}) and \cref{heat} (\cref{fig:err_heatold}). It is not surprising that they grow too at most linearly, since the model \cref{iterative} and the analysis given above are also applicable to the Chebyshev pseudospectral method. We can see that the rounding errors are comparable in these two methods and this should also be expected by the reasoning right above \cref{thm:bc}. 

\section{Computational cost} \label{sec:cost}
As pointed out in \cite{olv}, solving an almost banded system involves two steps: the QR factorization and the back substitution. They cost $\mO(m^2n)$ and $\mO(mn)$ respectively, where $n$ is the degrees of freedom and $m$ is the bandwidth of the almost banded matrix. 

The sparsity shown by \cref{fig:disc1} and \cref{fig:disc2} readily implies the same strategy for solving the resulting systems and, therefore, a computational cost of $\mO(n)$ too for both the discretization approaches regardless of the time marching scheme\footnote{For simplicity, we have omitted here the implied factor dependent of the bandwidth in the big-oh notation.}. However, when Approach 2 is employed we solve an upper triangular banded system (see \cref{fig:disc2}) for which only the back substitution is needed. Moreover, note that since the boundary condition \cref{Bucl} and the coefficients on the right-hand side of \cref{tulu} are independent of time, the QR factorization can be done once and for all at the beginning and at the subsequent steps only the back substitution is carried out.   

The $\mO(n)$ complexity is in stark contrast to the computational cost for solving \cref{tdpl} using the collocation-based pseudospectral method \cite{tre1}. If an explicit time marching scheme is used, the cost to calculate the derivatives on the right-hand side of \cref{tulu} is $\mO(n \log_2 n)$ with the aid of FFTs\footnote{In pseudospectral methods, variable coefficients are represented by diagonal matrices.}. If an implicit method is used, a dense system with no particular structure needs to be solved by a direct method such as LU factorization at a cost of $\mO(n^3)$ flops. Even though this cost is paid only once at the start of the time stepping and can be amortized over the subsequent steps, the cost of the backward substitution is still as high as $\mO(n^2)$ since the system is dense. Furthermore, for an adaptive implementation similar to the one introduced below in \cref{sec:adapt}, multiple or even a large number of LU decomposition may be needed, raising the cost significantly. These certify the great advantage of the ultraspherical spectral method in solving time-dependent PDEs.

\section{Adaptivity} \label{sec:adapt}
As time evolves, the solution to \cref{tdpl} may become spatially simpler or more complicated. It would be ideal if the method can take this into account and adapt the implementation for better efficiency but at the same time ensure that the degrees of freedom is large enough to guarantee an adequate resolution of the solution. This requires deciding a proper length of the solution vector at each time step. 

\begin{algorithm}
  \caption{Detection of a plateau \cite{aur}.}
  \label{alg:plateau}
  \begin{algorithmic}[1]
  \Procedure {plateau}{$u$, $tol$}
  \State{Step 1: Compute the normalized upper envelope of $u$.}
  \State{$envelope_j = \max_{j \leq k \leq n}{\lvert u_k\lvert}$}
  \If{$envelop_1 \ne 0$}
  \State{$envelope=envelope/envelope_1$}
  \EndIf
  \State{Step 2: Search for a plateau.}
  \For {$j \leftarrow 2, n$}
  \State{$j_2=round(1.25j+5)$}
  \State{$e_1 = envelope(j)$}
  \State{$e_2 = envelope(j_2)$}
  \State{$r = 3\left(1 - \log(e_1)/\log(tol)\right)$}
  \If{$\left(e_1 == 0 \mbox{ or } e_2/e_1 > r \right)$}
  \State{\textbf{return} $j$, $j_2$}
  \EndIf
  \EndFor
  \EndProcedure
  \end{algorithmic}
\end{algorithm}

Aurentz and Trefethen \cite{aur} propose an automated procedure in the context of function approximation for determining where to chop a Chebyshev series so that the truncated Chebyshev series is accurate and economical. The key of their chopping algorithm is the detection of a plateau, where the Chebyshev coefficients stay below a threshold and are sufficiently level. It is then based on this plateau that a chopping strategy is formulated. \cref{alg:plateau} summarizes the plateau detection part of their chopping algorithm. As we can see, when we approximate a given function by Chebyshev series the emergence of a plateau signals sufficient resolution, and a large portion of the plateau and all the trailing coefficients beyond the plateau are discarded for efficiency (not indicated in \cref{alg:plateau}). 

In the context of solving time-dependent PDEs, a plateau also serves as an indicator of adequate resolution. However, we only chop off the trailing coefficients beyond the plateau at each step. 

Suppose we are marching to the $(k+r)$th step using the information at $t^{k}, t^{k+1}, \ldots,$ $t^{k+r-1}$ by a multistep scheme or a Runge-Kutta method (for which $r = 1$). If there is no plateau in the computed solution $u^{k+r}$, we keep doubling the lengths of $u^{k}, u^{k+1}, \ldots, u^{k+r-1}$ by prolonging them with zeros and then re-calculate $u^{k+r}$ until a plateau emerges. This way, we come up with the following algorithm which allows adaptivity for the solution --- the solution vector is lengthened when an improved resolution may be effected or truncated when keeping some of the coefficients would not improve the resolution. 

\begin{algorithm}
\caption{Adaptive stepping from $t^{k}, t^{k+1}, \ldots, t^{k+r-1}$ to $t^{k+r}$.}
\label{alg:adapt}
\begin{algorithmic}[1]
\State{Stepping by a linear multistep or Runge-Kutta method to obtain the computed solution $u^{k+r}$. In case $u^k, u^{k+1}, \ldots, u^{k+r-1}$ are not of the same length, extend the shorter vectors to the length of the longest vector by prolonging them with zeros before stepping.}
\Procedure {adapt}{$u^k$, $u^{k+1}$, \ldots, $u^{k+r}$}
\State{L = \textsc{length}$(u^{k+r})$} \Comment{Function \textsc{length} returns the length of a vector.}
\State{Call \textsc{plateau}($u^{k+r}$, $tol$).}
\If {there is a plateau formed by $\{u^{k+r}_i\}_{i=j}^{j_2}$ }
\State{Drop $\{u^{k+r}_i\}_{i = j_2+1}^L$, use $u^{k+r} = \{u^{k+r}_i\}_{i = 0}^{j_2}$ for computation at future steps.}
\Else
\State{$u^k = [u^k, 0, \ldots, 0]$ (padding with zeros so that the lengths of $u^k$ is $2L$) for $k = 0, 1, \ldots, r-1$}
\State{Re-calculate $u^{k+r}$ by the same time marching scheme}
\State{Call \textsc{adapt}($u^k$, $u^{k+1}$, \ldots, $u^{k+r}$)}
\EndIf
\EndProcedure
\end{algorithmic}
\end{algorithm}

Note that the computed solution vectors fed into the calculation of the future 
steps are the ones with the plateau coefficients kept, i.e., only the trailing 
coefficients beyond the plateau are discarded. However, when a solution vector 
is no longer used for stepping, its plateau part can be safely dropped for 
saving storage, since keeping the plateau coefficients would not be of any help 
in improving the accuracy of the solution.

\begin{figure}[tbhp]
\centering
\subfloat[Lengths of the solution vectors: plateau coefficients included (solid) and plateau discarded (dotted).]{\label{fig:choppingvector}\includegraphics[scale=0.385]{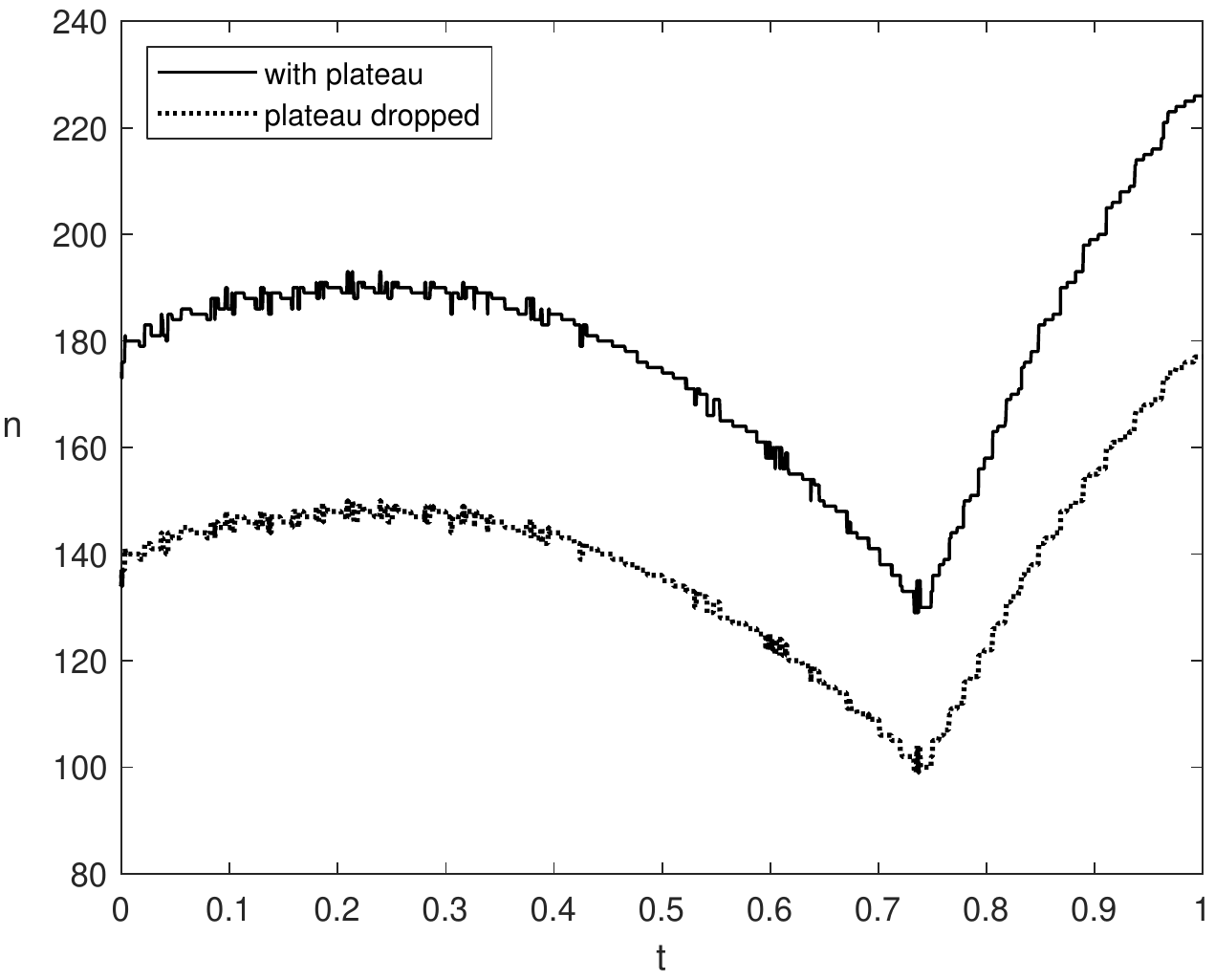}} \hspace{5mm}
\subfloat[Solution of \cref{variable}.]{\label{fig:choppingsolution}\includegraphics[scale=0.44]{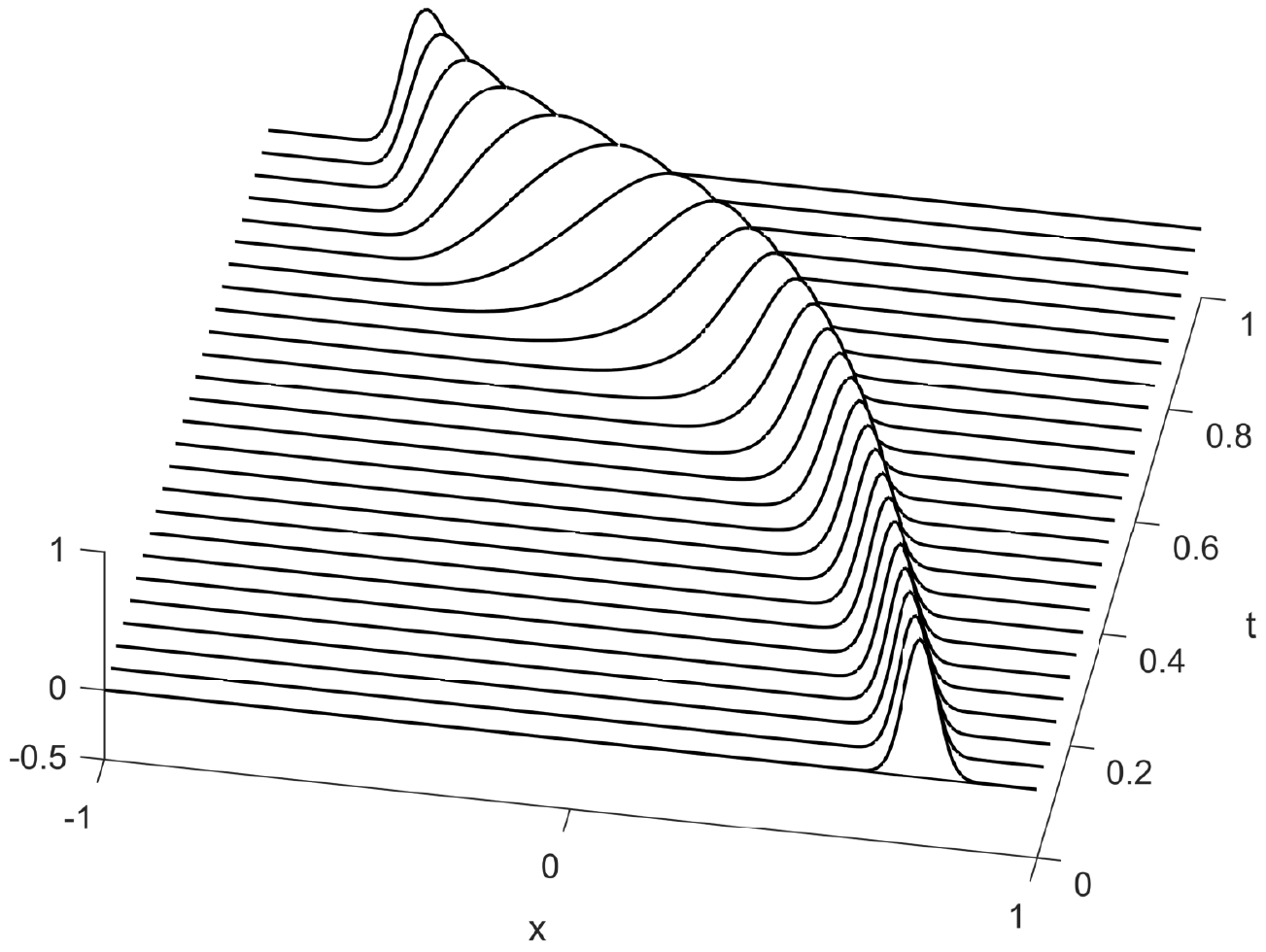}} 
\caption{Solving transport equation \cref{variable} of variable speed with adaptivity.}
\label{fig:chopping}
\end{figure}

To demonstrate how \cref{alg:adapt} works, we solve
\begin{align}
u_t = c(x) u_x  ~~\text{ s.t. } u(1, t) = 0, ~~ u(x, 0) = e^{-400(x-0.75)^2}, \label{variable}
\end{align}
where $c(x) = 3/5 + 3\sin^2(x-1)^2$ is a variable propagation speed depending on 
$x$ which results in a deformation of the left-travelling wave, as displayed in 
\cref{fig:choppingsolution}. The solid line in \cref{fig:choppingvector} shows 
the evolution of the length $n$ of the solution vector at each step, up to 
final time $t = 1$. This length includes the coefficients forming the plateau, 
whereas the dotted line shows the length if the plateau coefficients are 
discarded.

A noteworthy point is that the systems with different dimensions due to an 
adaptive implementation are not unrelated. Suppose that we solve with 
adaptivity and systems of dimensions $n_1 \times n_1$ and $n_2 \times n_2$ are 
solved by the QR factorization at two occasions with $n_2 > n_1$. Since the 
$n_2 \times n_2$ system is plainly an augment of the $n_1 \times n_1$ one by 
$n_2-n_1$ more rows and columns, we can simply cache the QR factorization for 
whichever system comes first to speed up the calculation for the other. Hence, 
for an adaptive implementation with systems of various sizes, the actual cost 
could be as little as doing the QR factorization once --- only for the system 
with the largest dimension.

Finally, we note that the \textsc{Chebfun} system \cite{dri}, particularly its 
PDE solver \textsc{pde15s}, offers a similar adaptivity in space. However, it 
is much more basic than the proposed one in that it does not have a mechanism 
for reducing the degrees of freedom when it is larger than it needs to be. 
Therefore, over-resolution may cause unnecessary drag in speed when the 
solution becomes spatially smoother. 

\section{Spatially periodic problems} \label{sec:periodic}
Up to this point, our discussion has been concentrated on spatially non-periodic problems. For \cref{ode} subject to periodic boundary conditions, one can simply follow the same framework of \cite{olv} but take $\left\{e^{\pm ikx}\right\}_{k = 0}^{\infty}$ as the basis functions, reproducing the tau-method \cite{ort}. With the Fourier basis, the $\lambda$-order differentiation operator remains sparse as 
\begin{align*}
\mD_{\lambda} = \displaystyle \diag\left(0 , i^{\lambda},(-i)^{\lambda},(2i)^{\lambda},(-2i)^{\lambda} ,\ldots , (2k)^{\lambda} ,(-2k)^{\lambda}, \ldots  \right) 
\end{align*}
and there is no more need for conversion operators $\mS_{\lambda}$, resulting 
in an even simpler implementation of the ultraspherical spectral method in the 
periodic case. However, for time-dependent PDEs with periodic boundary 
conditions, i.e., \cref{tdpl} with \cref{Bucl} replaced by periodic boundary 
conditions, time marching is not as easy as in the non-periodic case. For the 
simple cases where the right-hand side of \cref{tulu} is an odd-order spatial 
derivative of $u$, all the eigenvalues of the spatial discretization matrix 
reside on the imaginary axis for which only the schemes with a stability region 
enclosing the origin and its neighborhood along the imaginary axis are 
applicable. For example, for the one-dimensional transport equation \cref{wave1} with 
periodic boundary conditions, this immediate disqualifies all the explicit 
Runge-Kutta methods, the first two Adams-Bashforth methods, the Adams-Moulton 
methods of 2, 3, and 4 steps, and the BDF methods with more than 2 steps.

\section{Nonlinearity} \label{sec:nonlinearity} 
So far, the discussion has been concentrated on linear problems, which help 
simplify the analysis substantially. We now return to \cref{tdp} where $\mF$ 
also includes a nonlinear part as in \cref{LN}. In the remainder of this 
article, we slightly abuse the notation by assuming that the nonlinear operator 
$\mF$ takes in and returns Chebyshev and $C^{(\lambda)}$ coefficients 
respectively, instead of function values. This way, the input and the output of 
$\mF$ are consistent with those of the linear part $\mL$. For Approach 1, the 
fully discretized system reads
\begin{equation}
    \begin{split}
        &\begin{pmatrix}
        \mB\mP_{n}^{\top}\\ 
        \mP_{n-N}{\mS}_{N-1}\ldots \mS_0 \mP_{n}^{\top}
        \end{pmatrix} \mP_{n}\bu^{k+r} \\
        & = \begin{pmatrix}
        c \\ 
        \sum\limits_{j = 0}^{r - 1}(\beta_j h\mP_{n-N} \mF(t_{k+j}, \mP_{n}\bu^{k + 
        j}) - \alpha_j \mP_{n-N}{\mS}_{N-1}\ldots \mS_0 \mP_{n}^{\top} 
        \mP_{n}\bu^{k + j})
        \end{pmatrix}, \label{nl1ex}        
    \end{split}
\end{equation}
if an explicit multistep method ($\beta_r = 0$) is used. For an implicit 
multistep method ($\beta_r \neq 0$), we end up with the nonlinear equation 
\begin{equation}
    \begin{split}
        &\begin{pmatrix}
        \mB\mP_{n}^{\top}\mP_n \bu^{k + r}\\ 
        \mP_{n-N}{\mS}_{N-1}\ldots \mS_0 \mP_{n}^{\top}\mP_{n}\bu^{k+r} - h\beta_r\mP_{n-N} \mF(t_{k+r} , \mP_{n}\bu^{k+r})
        \end{pmatrix}  \\&= \begin{pmatrix}
        c \\ 
        \sum\limits_{j = 0}^{r - 1}(\beta_j h\mP_{n-N} \mF(t_{k+j}, \mP_{n}\bu^{k + j}) - \alpha_j \mP_{n-N}{\mS}_{N-1}\ldots \mS_0 \mP_{n}^{\top} \mP_{n}\bu^{k + j})
        \end{pmatrix}.
    \end{split} \label{nl1im}
\end{equation} 
The last two equations should be compared with \cref{disc1_lms}. If a 
Runge-Kutta method is used, \cref{disc1_rk} should be adapted to become
\begin{align}
\begin{pmatrix}
\mB\mP_{n}^{\top}\\
\mP_{n-N}{\mS}_{N-1}\ldots \mS_0 \mP_{n}^{\top}
\end{pmatrix}y_j &= 
\begin{pmatrix}
0 \\
h\mP_{n-N}\mF(t_k + \theta_jh, \mP_{n}\bu^{k} + \mu_jy_{j-1}) 
\end{pmatrix}, \label{nl1rk}
\end{align}

For Approach 2, explicit and implicit multistep methods leads to
\begin{equation}
    \begin{split}
        \mP_{n}&{\mS}_{N-1}\ldots \mS_0 \mP_{n}^{\top} \mP_{n}\bu^{k + r} \\
        & = h\sum_{j = 
        0}^{r - 1}\beta_j \mP_{n} \mF(t_{k+j} , \mP_{n} \bu^{k + j}) - 
        \sum_{j = 0}^{r-1}\alpha_j \mP_{n}{\mS}_{N-1}\ldots \mS_0 \mP_{n}^{\top} 
        \mP_{n}\bu^{k + j}
    \end{split} \label{nl2ex} 
\end{equation}
and 
\begin{equation}
    \begin{split}
        \mP_{n}&{\mS}_{N-1}\ldots \mS_0 \mP_{n}^{\top}\mP_{n}\bu^{k + r}- h\beta_r\mP_{n} 
        \mF(t_{k + r}, \mP_{n}\bu^{k + r}) \\
        & = h\sum_{j = 0}^{r - 
        1}\beta_j \mP_{n} \mF(t_{k+j} , \mP_{n} \bu^{k + j}) - \sum_{j = 
        0}^{r-1}\alpha_j \mP_{n}{\mS}_{N-1}\ldots \mS_0 \mP_{n}^{\top} \mP_{n}\bu^{k + 
        j}, 
    \end{split} \label{nl2im}
\end{equation}
respectively, where Runge-Kutta methods gives
\begin{align}
\mP_{n}{\mS}_{N-1}\ldots \mS_0 \mP_{n}^{\top} y_j = h \mP_{n}\mF(t_k + \theta_jh, \mP_{n}\bu^{k} + \mu_j y_{j-1}), \label{nl2rk}
\end{align}
which should be contrasted with \cref{tulu_rk_a2}. The nonlinear part $\mN(t, 
u)$ of $\mF(t, u)$ at specific $t$ and $u$ is usually evaluated by plugging in 
the value of $t$\footnote{In practice, $\mN$ is often independent of $t$, 
being a univariate function of $u$.}, sampling $\mN(u(x))$ at Chebyshev grids 
in $x$ of increasing size, calculating the Chebyshev coefficients by FFT until 
complete resolution, and converting them to $C^{(\lambda)}$ coefficients. The 
total cost is dominated by the few FFTs for the value-to-coefficient transform. 
Thus, the nominal complexities for solving \cref{nl1ex}, \cref{nl1rk}, 
\cref{nl2ex}, and \cref{nl2rk} are all $\mO(n \log_2n)$, where the linear 
complexity of system solving is prevailed over by the complexity of the 
evaluation of the nonlinear terms. Note that \cref{nl1im} and \cref{nl2im} are 
nonlinear equations of $\mP_{n}\bu^{k + r}$ and the cost of solution may be the 
greatest concern since the multiplication operators lose bandedness. However, 
it has been shown that fast solution to the nonlinear systems obtained from 
ultraspherical discretization can still be effected with $\mO(n \log_2 n)$ 
flops per iteration by an inexact Newton-GMRES method \cite{qin}. Since the 
solution at the previous time step can always serve as a good initial iterate 
for the next step, Newton's method usually skips the global stage and converges 
to machine precision in very few iterations. Thus, ultraspherical spectral 
method guarantees fast solution for virtually all the scenarios -- explicit and 
implicit schemes, linear and nonlinear equations. This is in marked contrast to 
solving time-dependent PDEs with the collocation-based pseudospectral method as 
the corresponding differentiation matrices are dense and much less structured.

The convergence of the solutions obtained by the two approaches in the 
nonlinear case is guaranteed if $\mF(t, u)$ satisfies Lipschitz conditions. 
This is met by virtually all the real-world problems.

To see how rounding errors accumulate, we replace the iterative model 
\cref{iterative} by
\begin{align*}
U^{k+1} = g(U^{k}),
\end{align*}
where $g$ is the nonlinear map corresponding to $\mF$. It can be shown that the modulus of the rounding error
\begin{align*}
\lnorm \Delta^{k+1}\rnorm  = \lnorm \udU^{k+1} - g(\udU^k)\rnorm  = 
\lnorm fl(\udg(\udU^k))- g(\udU^k)\rnorm \leq C_1 \eps,
\end{align*}
where $\udg$ denotes the floating point approximation to $g$ and $C_1 = (2 + 
\eps)\lnorm g(\udU^{k})\rnorm$. The accumulative error $E^{K+1}$ at $(K+1)$th 
step is bounded by
\begin{align*}
\lnorm E^{K+1}\rnorm = \lnorm \sum_{j = 
1}^{K+1}g^{K+1-j}\left(\Delta^{j}\right)\rnorm \nonumber \leq 
C_2(K+1)\eps_{mach},
\end{align*}
where $C_2 = (2 + \eps)\sup_{0\le r\le K}\lnorm g^{r}(\cdot)\rnorm \sup_{1\le 
j\le K+1}\lnorm g(\udU^{j-1})\rnorm$. This constant $C_2$ is, again, solely 
determined by the nonlinear map $g$. The conclusion that the rounding error, in the worst possible scenario, renders a linear growth is unchanged.

How the adaptivity described in \cref{sec:adapt} is implemented is not affected by the nonlinearity and, thus, stays the same. For nonlinear periodic problems, the evaluation of the nonlinear term $\mF$ is done with the Fourier coefficients, analogous to their Chebyshev counterpart in a straightforward manner.

The implementation of more advanced methods, such as the implicit-explicit differencing method, shares substantial similarities with those of the multistep and the Runge-Kutta methods. We choose to omit the discussion here.

\section{Exponential integrators} \label{sec:ei}
We could have closed this article at the end of last section. But the exponential integrators, also known as exponential time differencing, deserve a detailed discussion -- it is arguably the most powerful method for solving stiff ODE initial value problems. More importantly, the combination of the 
exponential integrators and the ultraspherical spectral method turns out to be extremely efficient, as we shall see below.

Consider the main equation \eqref{tuF}, that is, 
\begin{align*}
\mT u = \mL u + \mN(t, u).
\end{align*}
Suppose the linear operator $\mL$ is expressed as in \eqref{L}. To ensure that the coefficients produced on both sides are in the same space, we premultiply the right-hand side by $\mS_{0}^{-1} \ldots \mS_{N-1}^{-1}$ to obtain
\begin{align*}
\mT \bu &= \mS_{0}^{-1} \ldots \mS_{N-1}^{-1} \mL \bu + \mN(t, \bu).
\end{align*}
Following Approach 2 in \cref{sec:disc}, we ignore the boundary conditions momentarily and integrate the last equation on both sides from $t_k$ to $t_{k + 1}$ to have the variation-of-constant formula in terms of the ultraspherical spectral operators
\begin{equation}
    \begin{split}
        \bu(t_{k + 1}) =& e^{h \mS_{0}^{-1} \ldots \mS_{N-1}^{-1} \mL}\bu(t_k) 
        + e^{h \mS_{0}^{-1} \ldots \mS_{N-1}^{-1} \mL} \\
        \times & \int_{0}^{h}e^{\tau \mS_{0}^{-1} 
        \ldots \mS_{N-1}^{-1} \mL} \mN(t_k + \tau , \bu(t_k + \tau)) \, d\tau. 
    \end{split}\label{etdmain}
\end{equation} 

Different approximations to the integral in \eqref{etdmain} lead to various classes of exponential integrator \cite{hoc2}. If the integrand in \eqref{etdmain} is replaced by its polynomial interpolant at certain distinct points in $[t_k, t_{k+1}]$, we have the exponential multistep methods
\begin{align}
\bu^{k + 1} = e^{h \mS_{0}^{-1} \ldots \mS_{N-1}^{-1} \mL} \bu^{k} + h \sum_{j = 0}^{p - 1} \zeta_{j}(h \mS_{0}^{-1} \ldots \mS_{N-1}^{-1} \mL) \nabla^{j} \bv^{k}, \label{ems}
\end{align}
where $\bu^k=\bu(t_k)$, $\bv^k = \mN(t_k, \bu^k)$, and $\nabla^{j}\bv^{k}$ denotes the $j$th backward difference defined recursively by $\nabla^{0}\bv^{k} = \bv^{k}$ and $\nabla^{j+1}\bv^{k} = \nabla^{j}\bv^{k} - \nabla^{j}\bv^{k-1}$. The weights $\zeta_{j}$ can be calculated via the recurrence relation
\begin{align*}
\zeta_{0}(z) &= \varphi_{1}(z),\\
z\zeta_{j}(z) + 1 &= \sum_{i = 0}^{j-1}\frac{1}{j-i}\zeta_{i}(z), \nonumber
\end{align*}
where $\varphi_{1}(z) = \frac{e^z - 1}{z} $ is one of the so called $\varphi$-functions. These $\varphi$-functions can be generated from $\varphi_0(z) = e^{z}$ and the recurrence relation \cite{hoc1}
\begin{align*}
\varphi_{j + 1}(z) = \frac{\varphi_{j}(z) - \varphi_{j}(0)}{z}.
\end{align*}

Similarly, replacing the integrand in \eqref{etdmain} by its Taylor expansion at $t_k$ gives the exponential Runge-Kutta methods
\begin{subequations}
\begin{align}
\bu^{k + 1} = \exp(h \mS_{0}^{-1} \ldots \mS_{N-1}^{-1} \mL)\bu^k + h\sum_{i = 1}^{s}b_i(h \mS_{0}^{-1} \ldots \mS_{N-1}^{-1} \mL)\bv^{ki}, \\
\bu^{ki} = \exp(c_{i} h \mS_{0}^{-1} \ldots \mS_{N-1}^{-1} \mL)\bu^k + h\sum_{j = 1}^{s}a_{ij}(h \mS_{0}^{-1} \ldots \mS_{N-1}^{-1} \mL) \bv^{kj},
\end{align}\label{erk}%
\end{subequations}
where $\bu^{ki} = \bu(t_k + c_ih)$, $\bv^{ki} = \mN(t_k + c_ih, \bu^{ki})$. Like the Runge-Kutta methods, the weights $a_{ij}$ and $b_i$ satisfy $\sum_{j = 1}^s b_{j}(z) = \varphi_{1}(z)$ and $\sum_{j = 1}^{s} a_{ij}(z) = c_{i}\varphi_{1}(c_{i}z)$ for $i = 1, 2, \ldots, s$. For \eqref{erk} to be explicit, it is also required that $c_{1} = 0$ and $a_{ij}(z) = 0$ for $1 \le i \le j \le s$. 

Various exponential Runge-Kutta methods have been constructed and some of the most commonly used higher-order schemes are those proposed by Cox and Matthews \cite{cox}, Krogstad \cite{kro}, and Hochbruck and Ostermann \cite{hoc1}. For example, the method by Krogstad is given by the following Butcher tableau
\begin{equation*}
    \begin{matrix}
        \begin{array}{c|cccc}
            \setlength\tabcolsep{20pt}
            c_1 = 0\,& & & & \\[3pt]
            c_2 = \frac{1}{2}\,& a_{21} = \frac{1}{2}\varphi_{1,2}& & & \\[3pt] 
            c_3 = \frac{1}{2}\,& a_{31} = \frac{1}{2}\varphi_{1,3} - \varphi_{2,3}& a_{32} = \varphi_{2,3}& & \\[3pt]
            c_4 = 1\,& a_{41} = \varphi_{1,4} - 2\varphi_{2,4}\;& & \hspace{6pt} a_{43} = 2\varphi_{2,4}& \\[3pt]
            \hline \\[-10pt]
            \,&\,b_1 = \varphi_{1} - 3\varphi_{2} + 4\varphi_{3}\hspace{6pt}& b_2 = 2\varphi_{2} - 4\varphi_{3}& b_3 = b_2& \hspace{6pt} b_4 = -\varphi_{2} + 4\varphi_{3}
        \end{array}  &
        \begin{array}{l}
            \,\\[3pt]
            \,\\[3pt]
            \,\\[12pt]
            \hspace{-5.5pt},\\
            \,\\
        \end{array}
    \end{matrix}
\end{equation*}
where $\varphi_{i,j}(z) = \varphi_{i}(c_jz)$. 

To make the exponential multistep method \eqref{ems} and exponential Runge-Kutta method \eqref{erk} practical, we still need to truncate all the operators and infinite vectors to finite dimensions. This is done by replacing $h\mS_{0}^{-1} \ldots \mS_{N-1}^{-1} \mL$ by $G = h \mP_n\mS_{0}^{-1} \ldots \mS_{N-1}^{-1} \mL \mP_n^{\top}$ and only retaining the first $n$ components of $\bu^k$, $\bv^k$, and $\bv^{ki}$. For convenience, we denote by $u^k$, $v^k$, and $v^{ki}$ respectively the vectors formed by the first $n$ components of $\bu^k$, $\bv^k$, and $\bv^{ki}$.

The implementation of the exponential multistep and Runge-Kutta methods reviewed 
above boils down to the calculation of the product $\varphi_{j}(G)\xi$, where we 
use $\xi$ to denote any of $u^k$, $v^k$, and $v^{ki}$. Since evaluating 
$\varphi_{j}(G)$ directly usually suffers from large cancellation errors for 
$G$ of small magnitude, the evaluation of $\varphi_{j}(G)$ should be done via 
the Dunford-Taylor integral \cite{kas}. It is further shown that 
\begin{align*}
\varphi_{j}(z) = \frac{1}{2\pi i} \int_{\Gamma}\frac{e^s}{s^j} \frac{1}{s - z} \md s,
\end{align*}
where $\Gamma$ is a closed contour enclosing all the eigenvalues of $G$ \cite{sch}. Replacing $\Gamma$ by a $\theta$-parameterized Hankel contour $\phi(\theta)$, such as a Talbot's contour \cite{tre2}, leads to
\begin{align*}
\varphi_{j}(z) = \frac{1}{2\pi i} \int_{-\infty}^{+\infty}\frac{e^{\phi(\theta)}}{\phi(\theta)^j} \,\frac{1}{\phi(\theta) - z} \phi'(\theta)\md \theta.
\end{align*}
By truncating the integration interval to $[-\pi, \pi]$ and approximating the integral by $q$-point trapezoidal rule, we have a $q$-term sum-of-pole approximation of $\varphi_j(z)$
\begin{align}
r^{CI}_{(j)}(z) = \sum_{l = 1}^{q}\frac{w_l^{CI}}{z - z_l^{CI}}, \label{ci2}
\end{align}
where $w_l^{CI} = iq^{-1}e^{\phi_l}\phi'_l / \phi_l^{j}$, $\phi_l = z_l^{CI} = \phi(\theta_l)$, $\phi'_l = \phi'(\theta_l)$, and $\theta_l = \pi(2l-q-1)/(q-1)$ for $l = 1, 2, \ldots, q$.

One can also use the Carath\'eodory-Fej\'er approximation \cite{tre2,tre3} to obtain a near-best rational approximation to $\varphi_j(z)$
\begin{align}
r^{CF}_{(j)}(z) = \sum_{l = 1}^{q} \frac{w^{CF}_l}{z - z^{CF}_l}, \label{cf}
\end{align}
which is also in the sum-of-pole form as the one found by contour integral. The poles $z^{CF}_l$ and weights $w^{CF}_l$ of the CF approximation to $\varphi_j(z)$ usually differ for different $j$. 

Note that when \eqref{ci2} or \eqref{cf} are used, the calculation of $\varphi_j(G)\xi$ turns to solving linear systems $(G - z_{l}I)x_l = \xi$, or equivalently
\begin{align}
(hL - z_{l}{S}_{N-1} \ldots S_0)x_l = {S}_{N-1}\ldots S_0 \xi, \label{eisys}
\end{align}
for $l = 1, 2, \ldots, q$. What makes the exponential integrator even more powerful in the current context is the fact that \eqref{eisys} is a banded system as $L$ and ${S}_{N-1} \ldots S_0$ are both banded. When the poles and weights are known from pre-computation, the total cost $\mO(qn)$ for computing each of $\varphi_j(G)\xi$ is significantly less than if the collocation-based pseudospectral method were used, for which \eqref{eisys} is dense. Note that the convergence rate of $r^{CF}_{(j)}(z)$ is twice of that of $r^{CI}_{(j)}(z) $ and further speed-up for the CF method can be achieved by using common poles for all $\varphi_j$, whereas the contour-based method can compute the weights and poles cheaply. For comparisons of the contour-based and the CF methods, see \cite{tre2,sch}.

Here is a quick example of exponential integrating the Fisher equation
\begin{align*}
u_t = 0.001u_{xx} + u - u^2, \quad x \in [-1, 1],
\end{align*}
subject to the homogeneous Dirichlet boundary conditions and the initial 
condition $u(0, x) = (1 - \tanh(40x/\sqrt{6}))/4$ using the ultraspherical and 
the pseudospectral spectral methods. For both the methods, we choose $n = 512$ 
and integrate up to $t = 10$ with steps of size $1/n^2$, contrasted with the 
$\mO(1/n^4)$ restriction derived in \cref{sec:stability}. It takes $2.7$ 
seconds for the ultraspherical spectral method to finish the simulation, which 
is compared with $3.9$ seconds using the collocation-based pseudospectral 
method. Two methods have comparable accuracy of $\mO(10^{-9})$ in this 
experiment. The acceleration is more substantial when the degrees of freedom is 
greater.

\begin{remark}
Different exponential integrators vary by how the integral of the nonlinear 
term is approximated. For a linear problem, i.e., $\mN = 0$, all the 
exponential integrators coincide and give the same solution
\begin{align}
u^{k+1} = \varphi_0(G)u^k. \label{expm}
\end{align}
Since this solution is exact, the step size is unlimited\footnote{In practice, 
$\varphi_0$ can hardly be evaluated accurately for extremely large argument due 
to the conditioning.}, and the exponential integrators are also superb in 
solving linear problems. For example, exponential integrating the heat equation
\begin{align*}
u_t = 0.1u_{xx},~~ u(0 , 1) = u(0 , -1) = 0,~~ u(x , 0) = \sin\left(2\pi 
x\right)
\end{align*}
by the ultraspherical spectral method with $n = 32$ allows the step size to be 
as large as $0.1$. With this step size, the absolute error of the computed 
solution at $t = 10$ is about $1.1352e-14$.

\textsc{Chebfun} has an \texttt{expm} function for calculating operator 
exponentials, which overrides the \textsc{Matlab} function with the same name 
but working on matrices. It can exactly be used for evaluating $\varphi_0(G)$ 
in \cref{expm}. However, \textsc{Chebfun} does not offer any more functionality 
in exponential integration beyond the linear case. Additionally, the 
\textsc{Chebfun} \texttt{expm} explicitly forms the matrix that approximates 
$\varphi_0(G)$ before it is applied to the vector $u^k$, therefore the sparsity 
seen in \cref{eisys} is not taken advantage of and the storage cost becomes 
$\mO(n^2)$ instead of $\mO(n)$.

\end{remark}

\section{Conclusion and remarks} \label{sec:conclusion}
We have applied the ultraspherical spectral method to solving time-dependent 
PDEs by offering two approaches for discretization and have examined a few key 
aspects of the proposed method, including the stability of stepping, the error 
accumulation, and the computational cost, for both the linear and nonlinear 
cases. Careful comparison shows that the new method ties with the Chebyshev 
pseudospectral method in terms of stability and error and has a clear advantage 
in speed and adaptivity. 

So far, we have seen banded or almost-banded systems in two scenarios -- the 
implicit multi-step methods like the Adam-Moulton and BDF methods and the 
exponential integrator. Since the sparsity is a consequence of the employment of 
the ultraspherical spectral method for the spatial discretization, many more 
time marching schemes can also enjoy the fast linear algebra when used for 
solving time-dependent PDEs. More advanced examples include the spectral 
deferred correction method \cite{dut} for obtaining high accuracy solutions, the 
parareal method for time integration in parallel \cite{gan}, the symplectic 
integrator for Hamiltonian systems \cite{rut}, just to name a few. When 
stiffness requires the use of basic implicit methods as the underlying driving 
schemes, the method benefits from the resulting sparse linear systems.

The speed-up that we have seen could be even more conspicuous for 
problems in higher spatial dimensions since the degrees of freedom $n$ is 
squared or cubed.

One thing we have left out but worth mentioning is the handling of a second 
derivative in time. If high accuracy is not required, it is usually 
approximated by the simple leap frog formula. A more general approach is to 
reduce an equation with a second-order temporal derivative to a system of two 
equations with first-order derivatives in time. For instance, $u_{tt} = 
\mF(t, u(x, t))$ is reduced to
\begin{align*}
v_t &= \mF(t, u(x, t)),\\
u_t &= v(x, t),
\end{align*}
where the methods covered in the previous sections can be applied.

\begin{figure}[t!]
    \centering
    \subfloat[$t=0$]{\label{fig:heat2d0}\includegraphics[scale=0.35]{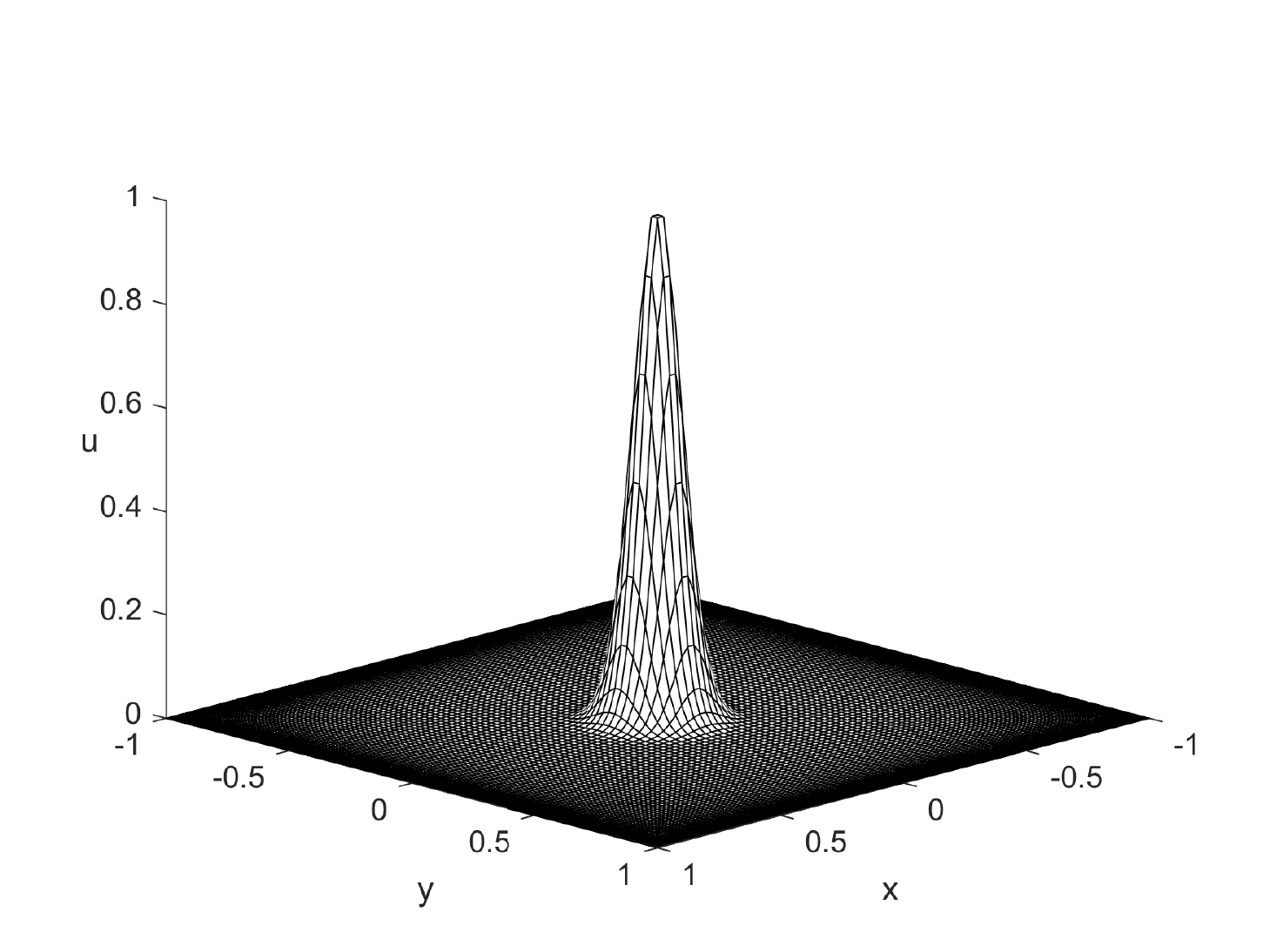}} \hspace{5mm}
    \subfloat[$t=1$]{\label{fig:heat2d1}\includegraphics[scale=0.35]{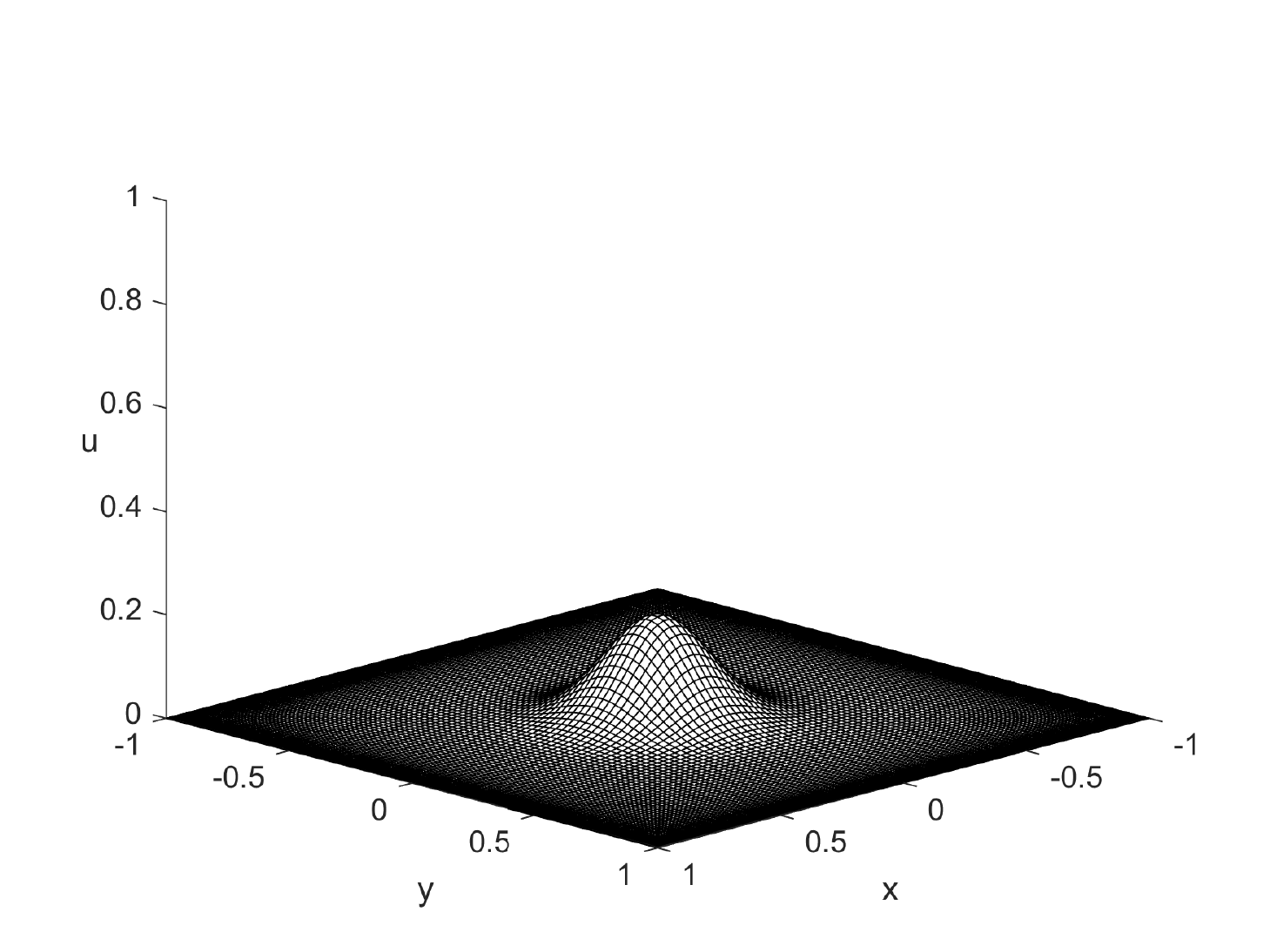}} 
    \caption{Solving a two dimensional heat equation by Approach 1 with $n_x = n_y = 128$ and the backward Euler with $h = 0.001$.}
    \label{fig:heat2d}
\end{figure}

Another possibility that is beyond the scope of this work is the extension of the ultraspherical spectral method to time-dependent problems in multiple spatial dimensions. The analysis may be more complicated and subtler than the present one, partly due to the boundary conditions. However, our initial numerical experiments show that the discretization approaches discussed in \cref{sec:disc} work well as expected. \cref{fig:heat2d1} displays the  solution at $t=1$ to the two-dimensional heat equation in a square domain subject to homogeneous boundary conditions
\begin{align*}
\begin{split}
   & u_t = 0.01\left(u_{xx} + u_{yy}\right), ~~~  (x, y) \in [-1 , 1] \times [-1 , 1],\\
   & \mbox{s.t.}~~~ u|_{\Gamma } = 0 ~~\text{ and }~~  u(x,0) = e^{-100(x^2 + y^2)},
    \end{split}
\end{align*}
where initial profile is shown in \cref{fig:heat2d0}.

As the ultraspherical spectral method has been widely accepted in the last decade, we believe the methods proposed in this article can serve as a natural companion of the ultraspherical spectral method for solving time-dependent problems and the analysis we have carried out can help understand and interpret the numerical results obtained from using these methods.

%
%

\bibliographystyle{spmpsci}      
\bibliography{bibliography}   

%
%

\end{document}